\documentclass[leqno, notitlepage]{amsart}

\usepackage{amssymb}
\usepackage[all, poly]{xy}

\DeclareMathOperator{\res}{res}
\DeclareMathOperator{\Hom}{Hom}

\DeclareMathOperator{\End}{End}

\DeclareMathOperator{\im}{im}

\DeclareMathOperator{\Id}{Id}

\DeclareMathOperator{\rank}{rank}

\newcommand{\sss}{\scriptscriptstyle}

\swapnumbers
\theoremstyle{plain} 
\newtheorem*{theorem}{Theorem}

\newtheorem*{lemma}{Lemma}
\newtheorem*{proposition}{Proposition}
\newtheorem{numtheorem}[subsection]{Theorem}

\theoremstyle{remark}
\newtheorem*{remark}{Remark}

\numberwithin{equation}{subsection}

\numberwithin{enumi}{subsection}
\newenvironment{alphenum}{
\begin{enumerate}
\setcounter{enumi}{\value{equation}}
}{
\setcounter{equation}{\value{enumi}}
\end{enumerate}}

\begin{document}

\title[Tensor product varieties]{
Tensor product varieties and crystals \\
GL case}

\author{ Anton Malkin }
\address{ Department of Mathematics, Yale University,
P.O. Box 208283, New Haven, CT 06520-8283 }
\email{anton.malkin@yale.edu}


\begin{abstract}
A geometric theory of tensor 
product for 
$\mathfrak{gl}_{\sss N}$-crystals is described. 
In particular,
the role of Spaltenstein varieties
in the tensor product is explained. 
As a corollary
a direct (non-combinatorial) proof
of the fact that the number of
irreducible components of a
Spaltenstein variety is equal 
to a Littlewood-Richardson
coefficient (i.e. certain
tensor product multiplicity)
is obtained.
\end{abstract}

\maketitle

\setcounter{section}{-1}

\setcounter{tocdepth}{1}
\tableofcontents

\section{Introduction}

\subsection{}

Given a partition $\lambda$ of an integer
$k$ into $N$ parts 
(zero parts allowed)
one can associate to it
various mathematical objects. Three of them
are of particular
interest to the representation theory:
a nilpotent orbit $O_{\lambda}$
in $\End_{\sss \mathbb{C}} \mathbb{C}^k$
(consisting of operators 
with the Jordan form given by $\lambda$), 
an irreducible polynomial representation 
$L (\lambda)$ of $GL_{\sss N}$ 
(with the highest weight $\lambda$),
and an irreducible representation 
$\rho_{\lambda}$ of the symmetric group 
$S_k$ over $\mathbb{C}$ 
(defined via the idempotent $e_\lambda$
in the group algebra $\mathbb{C}S_k$).
The goal of the geometric representation
theory in this context is to exploit 
the relation between nilpotent orbits
and representations of $S_k$ and
$GL_{\sss N}$.

\subsection{}

Let $t \in O_{\lambda}$ (i.e. $t$ is
a nilpotent operator in $\mathbb{C}^k$
with Jordan form $\lambda$), and let
$\mu^1 , \ldots , \mu^l$ be partitions
of integers $k^1 , \ldots , k^l$ into
$N$ parts each. Assume that
$k^1 + \ldots + k^l = k$, and consider
a variety 
$\mathfrak{S}_l
((\mu^1 , \ldots , \mu^l), \lambda)$
of all $l$-step partial
flags in $\mathbb{C}^k$ with dimensions
of the subfactors given by 
$k^1 , \ldots , k^l$ and such that
$t$ preserves each subspace in the flag, 
and when restricted to the subfactors
it defines  operators with the
Jordan forms given by
$\mu^1 , \ldots , \mu^l$. 
$\mathfrak{S}_l 
((\mu^1 , \ldots , \mu^l), \lambda)$
is called Spaltenstein variety. 
It is proven by Spaltenstein
(cf. \cite[II.5]{Spaltenstein1982}) that
$\mathfrak{S}_l
((\mu^1 , \ldots , \mu^l), \lambda)$
has pure dimension.
Let $\mathcal{S}_l
((\mu^1 , \ldots , \mu^l), \lambda)$
be the set of irreducible
components of $\mathfrak{S}_l 
((\mu^1 , \ldots , \mu^l), \lambda)$.
It follows from a theorem of Hall
\cite{Hall1959} (see also
\cite[Chapter II]{Macdonald1995})
that the cardinal of
$\mathcal{S}_l 
((\mu^1 , \ldots , \mu^l), \lambda)$
(i.e. the number of irreducible
components of
$\mathfrak{S}_l
((\mu^1 , \ldots , \mu^l), \lambda)$)
is equal to the multiplicity of
the representation $L_{\lambda}$
of $GL_m$ in the tensor product
$L_{\mu^1} \otimes \ldots \otimes 
L_{\mu^l}$:
\begin{equation}\label{IntroProduct}
| \mathcal{S}_l
((\mu^1 , \ldots , \mu^l), \lambda )| =
\dim_{\sss \mathbb{C}} 
\Hom_{\sss GL_{\sss N}}
(L_{\mu^1} \otimes \ldots \otimes 
L_{\mu^l} , \; L_{\lambda}) \; ,
\end{equation} 
or to the multiplicity of
$\rho_{\mu^1} \otimes \ldots \otimes 
\rho_{\mu^l}$ in the restriction of
$\rho_{\lambda}$ to
$S_{k^1} \times \ldots \times
S_{k^l}$:
\begin{equation}\label{IntroRestriction}
| \mathcal{S}_l
((\mu^1 , \ldots , \mu^l), \lambda )| =
\dim_{\sss \mathbb{C}} 
\Hom_{\sss S_{k^1} \times \ldots \times
S_{k^l}}
(\rho_{\mu^1} \otimes \ldots \otimes 
\rho_{\mu^l} , \; 
\res_{\sss S_k}^{\sss S_{k^1} \times \ldots 
\times S_{k^l}} \rho_{\lambda}) \; .
\end{equation} 
The right hand sides of 
\ref{IntroProduct} and 
\ref{IntroRestriction} are equal
due to Schur-Weyl duality, and
either equality can be derived
by combinatorial means (for example
from the Littlewood-Richardson Rule
-- cf. \cite[Chapter II]{Macdonald1995}).
However it is interesting to understand
the geometric meaning of 
\ref{IntroProduct} and 
\ref{IntroRestriction} (i.e to find
the role of the variety 
$\mathfrak{S}_l 
((\mu^1 , \ldots , \mu^l), \lambda )$ in
the tensor product for $GL_{\sss N}$
and in the restriction functor for 
symmetric groups).

The statement
\ref{IntroRestriction} was proven
geometrically by Borho and MacPherson
\cite{BorhoMacPherson} using
Springer's realization of representations
of symmetric groups. 

The goal of this paper is to develop a geometric
theory of tensor product for
polynomial representations of 
$GL_{\sss N}$, or, more exactly,
for $gl_{\sss N}$-crystals. 
In this geometric setup the variety $\mathfrak{S}_l 
((\mu^1 , \ldots , \mu^l), \lambda )$
is responsible for the multiplicities
in the tensor product decomposition.
In particular, \ref{IntroProduct} 
is an immediate corollary of the
main theorem of the paper
(Theorem \ref{GLnTheorem}).

\subsection{}

Let $t$ be again a nilpotent operator
in $\mathbb{C}^k$ with Jordan form
$\lambda$ (i. e. $t \in O_{\lambda}$).
Consider the variety 
$\mathfrak{M}_{\sss N} (\lambda)$
consisting of all $N$-step partial
flags
$(\{ 0 \} = F_0 \subset
F_1 \subset \ldots \subset
F_{\sss N}=\mathbb{C}^k)$
in $\mathbb{C}^k$ such that
$t F_k \subset F_{k-1}$ for
any $k = 1 , \ldots , N$. Note
that the dimensions of the subfactors
of the flag are not fixed. Each connected
component of the variety
$\mathfrak{M}_{\sss N} (\lambda)$
is of pure dimension.
Ginzburg \cite{Ginzburg1991}
constructed a representation of
$GL_{\sss N}$ in the top homology
of $\mathfrak{M}_{\sss N} (\lambda)$
(more exactly, in the sum of the
top Borel-Moore homology groups 
of connected components of
$\mathfrak{M}_{\sss N} (\lambda)$),
and proved that this representation is
isomorphic to $L (\lambda)$.
This paper deals with a weaker
structure -- that of
Kashiwara's 
$gl_{\sss N}$-crystal on the
set $\mathcal{M}_{\sss N} (\lambda)$
of irreducible components
of $\mathfrak{M}_{\sss N} (\lambda)$.
This crystal is isomorphic to the 
crystal $\mathcal{L} (\lambda )$
of the canonical basis of
$L (\lambda)$.

\subsection{}

The main idea of this paper
is to construct a certain variety 
$\mathfrak{T}_{\sss N} 
(\mu^1 , \ldots , \mu^l)$
such that its connected components
are of pure dimension and the set 
$\mathcal{T}_{\sss N} 
(\mu^1 , \ldots , \mu^l)$
of its irreducible components
can be equipped with a structure of
$gl_{\sss N}$-crystal in such a way that
one has two crystal isomorphisms:
\begin{equation}\label{IntroFirst}
\mathcal{T}_{\sss N} 
(\mu^1 , \ldots , \mu^l) \approx
\mathcal{M}_{\sss N} (\mu^1) 
\otimes \ldots \otimes 
\mathcal{M}_{\sss N} (\mu^l)
\end{equation} 
and
\begin{equation}\label{IntroSecond}
\mathcal{T}_{\sss N}
(\mu^1 , \ldots , \mu^l) \approx
\bigoplus_{\lambda} \mathcal{S}_l
((\mu^1 , \ldots , \mu^l), \lambda )
\otimes 
\mathcal{M}_{\sss N} (\lambda ) \; , 
\end{equation}
where the set $\mathcal{S}_l
((\mu^1 , \ldots , \mu^l), \lambda )$ is
considered as a trivial crystal,
and $\otimes$ (resp. $\oplus$)
denotes the tensor product
(resp. direct sum) of
crystals, which is equal
to the direct product
(resp. disjoint union)
on the level of sets.

The variety $\mathfrak{T}_{\sss N} 
(\mu^1 , \ldots , \mu^l)$ is
defined as the variety of all triples 
$(t, \mathbf{X}, \mathbf{F})$
consisting of a nilpotent operator 
$t$ in $\End_{\sss \mathbb{C}} 
\mathbb{C}^k$, an
$l$-step partial flag 
$\mathbf{X}$ in $\mathbb{C}^k$,
and an $N$-step partial flag 
$\mathbf{F}$ in $\mathbb{C}^k$,
such that subspaces of 
both flags are preserved
by $t$, and when $t$ is restricted
to the subfactors of $\mathbf{X}$
(resp. $\mathbf{F}$) it gives
operators with the Jordan forms
$\mu^1 , \ldots , \mu^l$
(resp. $0$ operators).
In particular the dimensions
of the subfactors of $\mathbf{X}$ 
are equal to $k^1, \ldots , k^l$, 
but the dimensions of the subfactors 
of $\mathbf{F}$ can vary.
The variety $\mathfrak{T}_{\sss N} 
(\mu^1 , \ldots , \mu^l)$ is called 
tensor product variety.

To describe the bijection 
\ref{IntroFirst} restrict the flag
$\mathbf{F}$ onto the subfactors
of the flag $\mathbf{X}$. 
In this way one obtains a map
from $\mathfrak{T}_{\sss N} 
(\mu^1 , \ldots , \mu^l)$ to
$\mathfrak{M}_{\sss N} (\mu^1)
\times \ldots \times 
\mathfrak{M}_{\sss N} (\mu^l)$,
which induces the bijection
\ref{IntroFirst}.

To obtain the bijection \ref{IntroSecond}
consider $\mathfrak{T}_{\sss N} 
(\mu^1 , \ldots , \mu^l)$ as a 
fibration over the nilpotent cone
in $\End_{\sss \mathbb{C}} 
\mathbb{C}^k$ (the projection map
takes $(t, \mathbf{X}, \mathbf{F})$
to $t$). The fiber of this 
fibration over a point 
$t \in O_{\lambda}$ is isomorphic to
the product of varieties
$\mathfrak{S}_l ((\mu^1 , \ldots ,
\mu^l), \lambda)$ and
$\mathfrak{M}_{\sss N} (\lambda)$.
Moreover the dimension of
the preimage of $O_{\lambda}$
in a connected 
component of $\mathfrak{T}_{\sss N} 
(\mu^1 , \ldots , \mu^l)$
does not depend on 
$\lambda$. Hence one have the
bijection \ref{IntroSecond}.

\subsection{}

The most important result of the paper
is that the crystal structures
on $\mathcal{T}_{\sss N} 
(\mu^1 , \ldots , \mu^l)$ induced
by the bijections 
\ref{IntroFirst} and \ref{IntroSecond}
are the same,
in other words that
the composite bijection
\begin{equation}\label{IntroTau}
\tau_{\sss N} : \quad
\mathcal{M}_{\sss N} (\mu^1) 
\otimes \ldots \otimes 
\mathcal{M}_{\sss N} (\mu^l)
\xrightarrow{\sim}
\bigoplus_{\lambda} \mathcal{S}_l
((\mu^1 , \ldots , \mu^l), \lambda )
\otimes 
\mathcal{M}_{\sss N} (\lambda )
\end{equation}
is a crystal isomorphism.
This is shown by reducing the problem to 
$GL_2$ case.

Actually in the main body of the paper only 
the case of two multiples (i.e. $l=2$)
is considered. However generalization
to arbitrary $l$ is straightforward.

\subsection{}

Certainly another proof of \ref{IntroProduct}
would not worth the effort. However
the understanding of its geometric
meaning paves the way for generalizations 
to Kac-Moody algebras other than 
$gl_{\sss N}$. In
particular in \cite{Malkin2000b}
the ADE case is considered, and 
certain ``multiplicity variety''
is defined such that the number
of its irreducible components is
equal to the multiplicity of 
a simple finite dimensional 
representation of a Lie algebra of ADE type
in the tensor product of $l$ simple
finite dimensional representations.
The multiplicity variety is
a generalization of the 
Spaltenstein variety in the context
of Nakajima's theory of quiver
varieties.

\subsection{}

The paper is organized as follows.
Section 1 contains a short review
of Kashiwara's theory of crystals.
In Section 2 the geometric tensor
product is described in  
$GL_2$ case. This 
demonstrates general constructions
on a simple example, but more important 
reason for including this section 
is that the generic $GL_{\sss N}$
case is reduced later to $N=2$.
In Section 3 the case of
$GL_{\sss N}$ is considered.
Both Sections 2 and 3 follow the
same pattern. First tensor product and
Spaltenstein varieties are defined,
then the bijections $\alpha$ and $\beta$
are described using certain 
``tensor product'' diagrams, and
finally it is proven that the
composite bijection $\tau$ is a crystal morphism
(in the case of $N=2$ by a direct verification
and for arbitrary $N$ via reduction to
$gl_2$-subalgebras inside $gl_{\sss N}$).

Most of the proofs consist of tedious
checks of dimensions of various
locally closed subsets of the varieties
involved. The only non-routine result is
\ref{GL2TTMMrank} which is the reason
for $\tau_{\sss N}$ to be a crystal
morphism.

\subsection{}

The tensor product variety in the special
case
$\mu^1 =(k^1, 0, 0, \dots)$, $\ldots$,
$\mu^l =(k^l, 0, 0, \dots)$,
was independently described by
H. Nakajima \cite[9.1]{Nakajima2001p}.

In the case when $k^1 = \ldots = k^l = 1$
(i.e. the case of the product of $l$ copies
of the fundamental representation of
$GL_{\sss N}$) the tensor product variety
is a certain Lagrangian subvariety of
the cotangent bundle to the variety
considered by Grojnowski and Lusztig
in \cite {GrojnowskiLusztig}. 

\subsection{}

Throughout the paper the following
conventions are used:
the ground field is $\mathbb{C}$;
``closed'', ``locally closed'', etc.
refer to the Zariski topology;
in the phrase
``locally trivial fibration'' 
``locally'' refers to the Zariski
topology, however trivialization
is analytic (not regular).

\section{Crystals}

Crystals were unearthed by Kashiwara 
\cite{Kashiwara1990,
Kashiwara1991, Kashiwara1994}.
An excellent survey of crystals as well
as some new results are given by 
Joseph in
\cite[Chapters 5, 6]{Joseph1995}.

\subsection{Definition of 
$\mathfrak{g}$-crystals}
Let $\mathfrak{g}$ be a reductive or 
a Kac-Moody Lie algebra, 
$I$ be the set of vertices of the
Dynkin graph of $\mathfrak{g}$,
$\mathcal{Q}_{\mathfrak{g}}$ 
(resp. 
$\mathcal{Q}^{\vee}_{\mathfrak{g}}$) 
be the weight (resp. coweight) 
lattice of $\mathfrak{g}$,
$\{ \alpha_i 
\in \mathcal{Q}_{\mathfrak{g}}
\}_{i \in I}$
(resp. 
$\{ \alpha^{\vee}_i \in 
\mathcal{Q}^{\vee}_{\mathfrak{g}} 
\}_{i \in I}$)
be the set of simple roots 
(resp. simple coroots),
$<,>$ be the natural
pairing between 
$\mathcal{Q}_{\mathfrak{g}}$
and $\mathcal{Q}^{\vee}_{\mathfrak{g}}$. 

A $\mathfrak{g}$-\emph{crystal} is a tuple 
$(\mathcal{A}, wt, 
\{ \varepsilon_i \}_{i \in I} ,
\{ \varphi_i \}_{i \in I} ,
\{ \Tilde{e}_i \}_{i \in I} ,
\{ \Tilde{f}_i \}_{i \in I} )$, where
\begin{itemize}
\item
$\mathcal{A}$ is a set,
\item
$wt$ is a map from $\mathcal{A}$ 
to $\mathcal{Q}_{\mathfrak{g}}$,
\item
$\varepsilon_i$ and $\varphi_i$ are maps from
$\mathcal{A}$ to $\mathbb{Z}$,
\item
$\Tilde{e}_i$ and $\Tilde{f}_i$ are maps 
from $\mathcal{A}$ to 
$\mathcal{A} \cup \{ 0 \}$. 
\end{itemize}
These data should satisfy the following axioms:
\begin{itemize}
\item
$wt (\Tilde{e}_i a) = wt (a) + \alpha_i$, 
\quad
$\varphi_i (\Tilde{e}_i a) = \varphi_i (a) + 1$, 
\quad
$\varepsilon_i (\Tilde{e}_i a) = 
\varepsilon_i (a) - 1$, 
\quad
for any $i \in I$ and $a \in \mathcal{A}$ 
such that $\Tilde{e}_i a \neq 0$,
\item
$wt (\Tilde{f}_i a) = wt (a) - \alpha_i$, 
\quad
$\varphi_i (\Tilde{f}_i a) = \varphi_i (a) - 1$, 
\quad
$\varepsilon_i (\Tilde{f}_i a) = 
\varepsilon_i (a) + 1$,
\quad
for any $i \in I$ and $a \in \mathcal{A}$ 
such that $\Tilde{f}_i a \neq 0$,
\item
$<\alpha^{\vee}_i , wt (a)> 
= \varphi_i (a) - \varepsilon_i (a)$
for any $i \in I$ and $a \in \mathcal{A}$,
\item
$\Tilde{f}_i a = b$ if and only if 
$\Tilde{e}_i b = a$, where
$i \in I$, and $a,b \in \mathcal{A}$.
\end{itemize}
The maps $\Tilde{e}_i$ and $\Tilde{f}_i$
are called Kashiwara's operators, and
the map $wt$ is called the weight function.

\begin{remark}
In a more general definition 
of crystals the maps $\varepsilon_i$ and 
$\varphi_i$ are allowed to have infinite values.
\end{remark}

A $\mathfrak{g}$-crystal 
$(\mathcal{A}, wt, 
\{ \varepsilon_i \}_{i \in I} ,
\{ \varphi_i \}_{i \in I} ,
\{ \Tilde{e}_i \}_{i \in I} ,
\{ \Tilde{f}_i \}_{i \in I} )$
is called \emph{trivial} if
\begin{equation}\nonumber
\begin{split}
wt (a) = 0 
&\text{ for any $a \in \mathcal{A}$,} \\  
\varepsilon_i (a) = \varphi_i (a) = 0 
&\text{ for any $i \in I$ and
$a \in \mathcal{A}$,} \\ 
\Tilde{e}_i a = \Tilde{f}_i a = 0 
&\text{ for any $i \in I$ and 
$a \in \mathcal{A}$.} \\ 
\end{split}
\end{equation}
Any set $\mathcal{A}$ can be equipped with the
trivial crystal structure as above.

A $\mathfrak{g}$-crystal 
$(\mathcal{A}, wt, 
\{ \varepsilon_i \}_{i \in I} ,
\{ \varphi_i \}_{i \in I} ,
\{ \Tilde{e}_i \}_{i \in I} ,
\{ \Tilde{f}_i \}_{i \in I} )$
is called \emph{normal} if
\begin{equation}\nonumber
\begin{split}
\varepsilon_i (a) &= \max 
\{ n \; | \; \Tilde{e}_i^n a \neq 0 \} \; , 
\\ \varphi_i (a) &= \max 
\{ n \; | \; \Tilde{f}_i^n a \neq 0 \} \; ,
\end{split}
\end{equation}
for any $i \in I$ and 
$a \in \mathcal{A}$. Thus in a normal 
$\mathfrak{g}$-crystal the maps 
$\varepsilon_i$ and $\varphi_i$
are uniquely determined by the action
of $\Tilde{e}_i$ and $\Tilde{f}_i$.
A trivial crystal is normal.

In the rest of the paper all 
$\mathfrak{g}$-crystals are assumed normal,
and thus the maps
$\varepsilon_i$ and $\varphi_i$ are usually
omitted.

By abuse of notation 
a $\mathfrak{g}$-crystal
$(\mathcal{A}, wt, 
\{ \varepsilon_i \}_{i \in I} ,
\{ \varphi_i \}_{i \in I} ,
\{ \Tilde{e}_i \}_{i \in I} ,
\{ \Tilde{f}_i \}_{i \in I} )$ is 
sometimes denoted simply by 
$\mathcal{A}$.

An \emph{isomorphism} of two 
$\mathfrak{g}$-crystals 
$\mathcal{A}$ and
$\mathcal{B}$
is a bijection between the sets 
$\mathcal{A}$ and $\mathcal{B}$ 
commuting with the action of
the operators $\Tilde{e}_i$ and
$\Tilde{f}_i$, and the functions 
$wt$, $\varepsilon_i$, and $\varphi_i$. 

The \emph{direct sum} 
$\mathcal{A} \oplus \mathcal{B}$
of two $\mathfrak{g}$-crystals 
$\mathcal{A}$ and $\mathcal{B}$
is their disjoint union as sets
with the maps $\Tilde{e}_i$,
$\Tilde{f}_i$, $wt$,
$\varepsilon_i$, and $\varphi_i$ 
acting on each
component of the union separately. 

\subsection{Tensor product of 
$\mathfrak{g}$-crystals}
The \emph{tensor product} 
$\mathcal{A} \otimes \mathcal{B}$
of two $\mathfrak{g}$-crystals 
$\mathcal{A}$ and $\mathcal{B}$
is their direct product as sets
equipped with the 
following crystal structure:
\begin{equation}
\label{DefinitionOfCrystalProduct}
\begin{split}
wt ((a,b)) &= wt (a) + wt (b)  \; ,
\\
\varepsilon_i ((a,b)) &=
\max \{ \varepsilon_i (a), \varepsilon_i (a)
+ \varepsilon_i (b) - \varphi_i (a) \} \; ,
\\
\varphi_i ((a,b)) &=
\max \{ \varphi_i (b), \varphi_i (a)
+ \varphi_i (b)- \varepsilon_i (b) \} \; ,
\\
\Tilde{e}_i ((a,b)) &=
\begin{cases}
(\Tilde{e}_i a \; , \; b)
&\text{ if } 
\varphi_i (a) \geq \varepsilon_i (b) \; ,\\
(a \; , \; \Tilde{e}_i b)
&\text{ if } 
\varphi_i (a) < \varepsilon_i (b) \; ,
\end{cases} 
\\
\Tilde{f}_i ((a,b)) &=
\begin{cases}
(\Tilde{f}_i a \; , \;  b)
&\text{ if } 
\varphi_i (a) > \varepsilon_i (b) \; ,\\
(a \; , \; \Tilde{f}_i b) 
&\text{ if } 
\varphi_i (a) \leq \varepsilon_i (b) \; .
\end{cases}
\end{split}
\end{equation}
Here $(a,0)$, $(0,b)$, and $(0,0)$ are
identified with $0$.
One can check that the set 
$\mathcal{A} \times \mathcal{B}$ with the
above structure satisfies all the axioms
of a (normal) crystal, and
that the tensor product of crystals
is associative. 

\begin{remark}
Tensor product of crystals 
is not commutative.
\end{remark}

\begin{remark}
The tensor product 
(in any order) of a crystal 
$\mathcal{A}$ with a trivial crystal 
$\mathcal{B}$ is isomorphic to the direct sum 
$\bigoplus_{b \in \mathcal{B}} \mathcal{A}_b$
where each $\mathcal{A}_b$ is isomorphic to 
$\mathcal{A}$.
\end{remark}

\subsection{Highest weight crystals and 
closed families}
A crystal $\mathcal{A}$ is a 
\emph{highest weight}
crystal with highest weight 
$\lambda \in \mathcal{Q}_{\mathfrak{g}}$ if
there exists an element 
$a_{\lambda} \in \mathcal{A}$ such that:
\begin{itemize}
\item 
$wt ( a_{\lambda} ) = \lambda$ ,
\item
$\Tilde{e}_i a_{\lambda} = 0$ 
for any $i \in I$,
\item
any element of $\mathcal{A}$ can be obtained
from $a_{\lambda}$ by successive applications
of the operators $\Tilde{f}_i$.
\end{itemize}

Consider a family of highest weight 
normal crystals
$\{ \mathcal{A} (\lambda ) \}_{\lambda 
\in \mathcal{J}}$ labeled by a set
$\mathcal{J} \subset 
\mathcal{Q}_{\mathfrak{g}}$ (the highest
weight of $\mathcal{A} ( \lambda )$ is
$\lambda$). The family
$\{ \mathcal{A} (\lambda ) \}_{\lambda 
\in \mathcal{J}}$ is called 
\emph{strictly closed} 
(with respect to tensor products)
if the tensor product of any two members of
the family is isomorphic to a direct sum 
of members of the family:
\begin{equation}\nonumber
\mathcal{A} ( \mu^1 ) \otimes
\mathcal{A} ( \mu^2 ) = 
\bigoplus_{\lambda \in \mathcal{J}}
\mathcal{U} ((\mu^1 , \mu^2 ), \lambda )
\otimes \mathcal{A} ( \lambda ) \; ,
\end{equation}
where 
$\mathcal{U} ((\mu^1 , \mu^2 ), \lambda )$
is a set with the trivial crystal structure.
A family
$\{ \mathcal{A} (\lambda ) \}_{\lambda 
\in \mathcal{J}}$ is called 
\emph{closed} 
if the tensor product 
$\mathcal{A} (\mu^1) \otimes
\mathcal{A} (\mu^2)$ 
of any two members of the family
contains $\mathcal{A}_{\mu^1+\mu^2}$
as a direct summand. Any strictly closed
family is closed. The converse is also true,
as a corollary of Theorem 
\ref{JosephTheorem}.

Let $\mathcal{Q}^{+}_{\mathfrak{g}}
\subset \mathcal{Q}_{\mathfrak{g}}$ be
the set of highest weights of integrable 
highest weights modules of $\mathfrak{g}$
(in the reductive case a module is called
integrable if it is derived from a 
polynomial representation of the corresponding
connected simply connected reductive group, 
in the Kac-Moody
case the highest weight should be a 
positive linear combination of the fundamental
weights). The original motivation for the 
introduction  of crystals was the discovery 
by Kashiwara \cite{Kashiwara1991} 
and Lusztig  \cite{Lusztig1991a}
of canonical (or crystal) bases in 
integrable highest 
weight modules of a (quantum) Kac-Moody algebra. 
These bases have many favorable properties, one
of which is that 
as sets they are equipped with a crystal 
structure. In other words to each 
irreducible integrable highest weight module
$L (\lambda )$ corresponds a normal crystal
$\mathcal{L} (\lambda )$
(crystal of the canonical basis). In this 
way one obtains a strictly closed
family of crystals
$\{ \mathcal{L} (\lambda ) \}_{\lambda 
\in \mathcal{Q}^+_{\mathfrak{g}}}$, 
satisfying the following two properties:
\begin{itemize}
\item
the cardinality of $\mathcal{L} (\lambda )$
is equal to the dimension of $L (\lambda)$;
\item
one has the following tensor product 
decompositions for 
$\mathfrak{g}$-modules and
$\mathfrak{g}$-crystals:
\begin{equation}\nonumber
\begin{split}
L ( \mu^1 ) \otimes
L ( \mu^2 ) &= 
\bigoplus_{\lambda 
\in \mathcal{Q}^+_{\mathfrak{g}}}
C ((\mu^1 , \mu^2 ), \lambda )
\otimes L ( \lambda ) \; ,
\\
\mathcal{L} ( \mu^1 ) \otimes
\mathcal{L} ( \mu^2 ) &= 
\bigoplus_{\lambda 
\in \mathcal{Q}^+_{\mathfrak{g}}}
\mathcal{C} ((\mu^1 , \mu^2 ), \lambda )
\otimes \mathcal{L} ( \lambda ) \; ,
\end{split}
\end{equation}
where the cardinality of the set
(trivial crystal)
$\mathcal{C} 
((\mu^1 , \mu^2 ), \lambda )$
is equal to the dimension of the linear space
(trivial $\mathfrak{g}$-module)
$C ((\mu^1 , \mu^2 ), \lambda )$.
\end{itemize}

\noindent
The aim of this paper is to 
construct another strictly closed 
family of $\mathfrak{g}$-crystals 
$\{ \mathcal{M} (\lambda ) \}_{\lambda 
\in \mathcal{Q}^+_{\mathfrak{g}}}$
(for $\mathfrak{g}$ being $gl_{\sss N}$
or a symmetric Kac-Moody algebra), using
geometry associated to $\mathfrak{g}$. The
following crucial theorem ensures that this
family is isomorphic to the family
$\{ \mathcal{L} (\lambda ) \}_{\lambda 
\in \mathcal{Q}^+_{\mathfrak{g}}}$
of crystals of canonical bases 
(two families of crystals labeled
by the same index set 
are called isomorphic
if the corresponding members of the
families are isomorphic as crystals). 
\begin{numtheorem}\label{JosephTheorem}
There exists a unique 
(up to an isomorphism) closed family of
$\mathfrak{g}$-crystals labeled
by $\mathcal{Q}^+_{\mathfrak{g}}$. 
\end{numtheorem}
\begin{proof}
For the proof in Kac-Moody case see
\cite[Proposition 6.4.21]{Joseph1995}. 
The statement for a reductive $\mathfrak{g}$ 
easily follows from the statement for
the factor of $\mathfrak{g}$ by its center.
\end{proof} 

\section{Grassmann varieties and $gl_2$-crystals}

This section contains a geometric description
of a closed family of $gl_2$-crystals.

\subsection{Notation}
From now on the weight lattice 
$\mathcal{Q}_{gl_2}$ of $gl_2$ is denoted 
simply by $\mathcal{Q}_2$ and identified
with $\mathbb{Z} \oplus \mathbb{Z}$.
In particular, $\mathcal{Q}_2^+ =
\mathbb{Z}_{\geq 0} \oplus \mathbb{Z}_{\geq 0}
\subset \mathbb{Z} \oplus \mathbb{Z}$.

Given $n,m \in \mathbb{Z}$, $Gr^n_m$ 
denotes the Grassmann variety of all $m$-dimensional 
subspaces of $\mathbb{C}^{n}$. If
$m < 0$ or $n < m$ then $Gr^n_m$ is empty.
Otherwise it is a smooth connected 
variety of dimension $m(n-m)$. 
Let $Gr (n) =\bigsqcup_{0 \leq m \leq n}
Gr_m^n$ be the variety of all subspaces
of $\mathbb{C}^n$.

\subsection{Geometric $gl_2$-crystals}
\label{GL2DefinitionOfM}
Given $w \in \mathbb{Z}_{\geq 0}$ consider 
the following diagram of varieties 
(cf. \cite{Ginzburg1991}):
\begin{equation}\label{GL2MpiM}
\xymatrix{
\mathfrak{M}_2 (w)
\ar[d]^{\pi_2} \\
\; \; \mathfrak{N}_2 (w) \; , }
\end{equation}
where

$\mathfrak{N}_2 (w) = 
\{ t \in \End (\mathbb{C}^w) 
\; | \; t^2 = 0 \}$,

$\mathfrak{M}_2 (w) = 
\{ (t, F) \; | \; t \in \mathfrak{N}_2 (w), \; 
F \in Gr(w), \;
\im t \subset F \subset \ker t \}$, in other
words, $\mathfrak{M}_2 (w)$ is a variety of
pairs $(t, F)$ consisting of an operator $t$ in
$\mathbb{C}^w$ and a subspace 
$F \subset \mathbb{C}^w$, such that 
$\im t \subset F \subset \ker t$,

$\pi_2 ((t, F)) = t$.

\noindent
The variety $\mathfrak{M}_2 (w)$ is 
a disjoint union of connected components:
\begin{equation}\nonumber
\mathfrak{M}_2 (w) = 
\bigsqcup_{v \in \mathbb{Z}_{\geq 0}}
\mathfrak{M}_2 (v,w)  \; ,
\end{equation}
where
\begin{equation}\nonumber
\mathfrak{M}_2 (v,w) = 
\{ (t,F) \in \mathfrak{M}_2 (w) \; | \; 
\dim F = v \}  \; .
\end{equation}
The variety $\mathfrak{N}_2 ( w )$ 
can be stratified as follows
\begin{equation}\nonumber
\mathfrak{N}_2 ( w ) =
\bigcup_{r \in \mathbb{Z}_{\geq 0}}
O_{( w-r,r )}  \; ,
\end{equation}
where
\begin{equation}\nonumber
O_{( w-r,r )} = 
\{ t \in \mathfrak{N}_2 
( w ) \; | \; \rank t = r \} \; .
\end{equation}
The weird way of writing index 
$\lambda$ of $O_{\lambda}$ will be
explained in subsection 
\ref{Spaltenstein}.
Each stratum $O_{\lambda}$ is a single 
$GL ( w )$-orbit in 
$\mathfrak{N}_2 ( w )$ and
there are finitely many non-empty strata.
The dimension of $O_{( w-r,r )}$ is given
by 
\begin{equation}\label{GL2DimO1}
\dim O_{(w-r,r)} =
2r (w-r) \; .
\end{equation}

Choose
$t \in O_{(w-r,r)}$ and consider the variety
$\mathfrak{M}_2 (w, r) = \pi_2^{-1} (t)$. 
This variety does not depend (up to an
isomorphism) on the choice of 
$t \in O_{(w-r,r)}$. However
if a specific $t$ is used then the notation
is $\mathfrak{M}_2 (w, t)$ instead of
$\mathfrak{M}_2 (w, r)$.

The variety $\mathfrak{M}_2 (w, r)$ is 
a disjoint union of connected components:
\begin{equation}\nonumber
\mathfrak{M}_2 (w, r) = 
\bigsqcup_{v \in \mathbb{Z}_{\geq 0}}
\mathfrak{M}_2 (v,w,r)  \; ,
\end{equation}
where
\begin{equation}\nonumber
\mathfrak{M}_2 (v,w,r) =
\mathfrak{M}_2 (w,r)
\cap
\mathfrak{M}_2 (v,w)  \; .
\end{equation}
The component $\mathfrak{M}_2 (v,w,r)$ is 
isomorphic 
to $Gr^{w-2r}_{v-r}$ (in particular it is
nonempty only for a finite number 
of values of $v$). Thus
$\mathfrak{M}_2 (v,w,r)$ is
a smooth connected variety of dimension
\begin{equation}\label{GL2DimM}
\dim \mathfrak{M}_2 (v,w,r) =
(w-r-v)(v-r) \; .
\end{equation}
Let $\mathcal{M}_2 (w,r) =
\{ \mathfrak{M}_2 (v,w,r) \}_{v = r}^{w-r}$
be the set of irreducible components 
of $\mathfrak{M}_2 (w,r)$. 
Endow this set with the
following structure of a
$gl_2$-crystal:
\begin{equation}
\label{GL2CrystalDefinition}
\begin{split}
wt (\mathfrak{M}_2 (v,w,r)) &= (v,w-v) \; ,
\\
\varepsilon (\mathfrak{M}_2 (v,w,r)) &= w-r-v \; ,
\\
\varphi (\mathfrak{M}_2 (v,w,r)) 
&= v-r \; ,
\\ 
\tilde{e} (\mathfrak{M}_2 (v,w,r)) &= 
\begin{cases} \mathfrak{M}_2 (v+1,w,r)
&\text{ if $v < w-r$} \; , \\
0 &\text{ if $v \geq w - r$} \; ,
\end{cases}
\\
\tilde{f} (\mathfrak{M}_2 (v,w,r)) &= 
\begin{cases} \mathfrak{M}_2 (v-1,w,r)
&\text{ if $v > r$} \; , \\
0 &\text{ if $v \leq r$} \; .
\end{cases}
\end{split}
\end{equation}
Here indexes of $\varepsilon$, $\varphi$,
$\Tilde{e}$, and $\Tilde{f}$ are omitted
because $gl_2$ has only one root.
The set $\mathcal{M}_2 (w,r)$
equipped with this structure is a highest
weight normal $gl_2$-crystal with the highest
weight $(w - r,r)$. In other words, it is 
isomorphic (as a crystal) to 
$\mathcal{L} ((w-r,r))$. 

The above is a trivial example of a geometric
construction of crystals due to Lusztig and 
Nakajima (cf. \cite{Lusztig1991a, Nakajima1998,
KashiwaraSaito}). Next subsection provides 
a geometric description 
of the crystal tensor product on the family
$\{ \mathcal{M}_2 (w,r ) 
\}_{w,r \in \mathbb{Z}_{\geq 0}}$.

\subsection{The tensor product variety}
\label{GL2DefinitionOfT}
Let $w^1, w^2, r^1, r^2$ be non-negative 
integers, and let 
$\mathfrak{T}_2 (w^1, r^1, w^2, r^2)$
be a variety of triples
$(t,X,F)$, where

$t \in \mathfrak{N}_2 ( w^1 + w^2 )$;

$F$ is a subspace of $\mathbb{C}^{w^1 + w^2}$,
such that $\im t \subset F \subset \ker t$
(i.e. $(t,F) \in \mathfrak{M}_2 (w^1 + w^2)$);

$X$ is a subspace in $\mathbb{C}^{w^1 + w^2}$
such that $\dim X = w^1$, $t X \subset X$, 
$\rank (t |_{X}) = r^1$, and
$\rank (t |_{(\mathbb{C}^{w^1 +w^2}/X)})=r^2$.

The variety 
$\mathfrak{T}_2 (w^1, r^1, w^2, r^2)$
is called \emph{tensor product variety}
for its role in the description
of the $gl_2$-crystal
$\mathcal{M}_2 (w^1, r^1)
\otimes
\mathcal{M}_2 (w^2, r^2)$.
One has the following decomposition of
$\mathfrak{T}_2 (w^1, r^1, w^2, r^2)$ 
into a disjoint union of varieties:
\begin{equation}\nonumber
\mathfrak{T}_2 (w^1, r^1, w^2, r^2) =
\bigsqcup_{v \in \mathbb{Z}_{\geq 0}} 
\mathfrak{T}_2 (v, w^1, r^1, w^2, r^2) 
\; ,
\end{equation}
where
\begin{equation}\nonumber
\mathfrak{T}_2 (v, w^1, r^1, w^2, r^2) =
\{ (t, X, F) \in 
\mathfrak{T}_2 (w^1, r^1, w^2, r^2)
\; | \; 
\dim F = v \} \; .
\end{equation}
Moreover 
\begin{equation}\nonumber
\mathfrak{T}_2 (w^1, r^1, w^2, r^2) =
\bigcup_{v^1, v^2 \in \mathbb{Z}_{\geq 0}}
\mathfrak{T}_2 (v^1, w^1, r^1, v^2, w^2, r^2) 
\; ,
\end{equation}
where
\begin{multline}\nonumber
\mathfrak{T}_2 (v^1, w^1, r^1, v^2, w^2, r^2) 
= \\ =
\{ (t, X, F) \in 
\mathfrak{T}_2 (w^1, r^1, w^2, r^2)
\; | \; 
\dim F \cap X = v^1 , \;
\dim F/(F\cap X) = v^2 \} \; .
\end{multline}
Each subset 
$\mathfrak{T}_2 (v^1, w^1, r^1, v^2, w^2, r^2)$
is locally closed and one has
\begin{equation}\nonumber
\mathfrak{T}_2 (v, w^1, r^1, w^2, r^2) =
\bigcup_{\substack{
v^1, v^2 \in \mathbb{Z}_{\geq 0} \\
v^1 + v^2 = v}}
\mathfrak{T}_2 (v^1, w^1, r^1, v^2, w^2, r^2) 
\; .
\end{equation}

In the next several subsections the variety
$\mathfrak{T}_2 (w^1, r^1, w^2, r^2)$ 
is used to describe the tensor product of 
$gl_2$-crystals
$\mathcal{M}_2 (w^1, r^1)$ and 
$\mathcal{M}_2 (w^2, r^2)$.
The set of irreducible components of
$\mathfrak{T}_2 (w^1, r^1, w^2, r^2)$ 
appears to be in bijection with the set 
$\mathcal{M}_2 (w^1, r^1) \times 
\mathcal{M}_2 (w^2, r^2)$, and two
ways of labeling the 
irreducible components 
imply the decomposition of the tensor
product as a direct sum.

\subsection{The tensor product variety 
as a fibered product}\label{GL2FibredProduct}
The variety $\mathfrak{N}_2 (w^1 + w^2)$ 
can be stratified by $GL (w^1 + w^2)$-orbits
(cf. subsection \ref{GL2DefinitionOfM}):
\begin{equation}\nonumber
\mathfrak{N}_2 (w^1 + w^2) =
\bigcup_{r \in \mathbb{Z}_{\geq 0}}
O_{(w^1+w^2-r,r)}  \; ,
\end{equation}
where
\begin{equation}\nonumber
O_{(w^1+w^2-r,r)} = 
\{ t \in \mathfrak{N}_2 
(w^1 + w^2) \; | \; 
\rank t = r \} \; . 
\end{equation}
The dimension of $O_{(w^1+w^2-r,r)}$ 
is given by \eqref{GL2DimO1}: 
\begin{equation}\label{GL2DimO2}
\dim O_{(w^1+w^2-r,r)} =
2r (w^1 + w^2 -r) \; .
\end{equation}
Consider a diagram
\begin{equation}\nonumber
\xymatrix{
\mathfrak{T}_2 (v, w^1, r^1, w^2, r^2) 
\ar[d]^{\delta_2} \\
\mathfrak{N}_2 (w^1 + w^2) \; ,
}
\end{equation}
where $\delta_2 ( (t, X, F))=t$.
The map $\delta_2$ restricted to a stratum
$O_{(w^1+w^2-r,r)}$ is a locally trivial
fibration with a fiber isomorphic to
a direct product 
\begin{equation}\nonumber
\mathfrak{M}_2 (v, w^1 + w^2,r) \times
\mathfrak{S}_2 (((w^1, r^1), (w^2, r^2)), 
(w^1 + w^2 ,r)) \; ,
\end{equation}
where the variety
$\mathfrak{S}_2 (((w^1, r^1), (w^2, r^2)), 
(w^1+ w^2 ,r))$
can be described as follows. Fix 
$t \in O_{(w^1 + w^2 -r, r)}$. Then
$\mathfrak{S}_2 (((w^1, r^1), (w^2, r^2)), 
(w^1 + w^2,r))$
is a subvariety of the Grassmannian 
$Gr_{w^1}^{w^1 + w^2}$
of $w^1$-dimensional subspaces of 
$\mathbb{C}^{w^1+w^2}$ given by
\begin{multline}\nonumber
\mathfrak{S}_2 (((w^1, r^1), (w^2, r^2)), 
(w^1 + w^2,r)) 
= \\
\{ X \subset \mathbb{C}^{w^1 + w^2} \; | \;
\dim X = w^1 , \;
t X \subset X  , \;
\rank ( t |_X )  = r^1 , \;
\rank ( t |_{\mathbb{C}^{w^1 + w^2}/X} ) 
= r^2    \} 
\end{multline}
The variety 
$\mathfrak{S}_2 (((w^1, r^1), (w^2, r^2)), 
(w^1 + w^2,r))$ is an 
example of a Spaltenstein variety
(cf. subsection \ref{Spaltenstein}). 
It does not depend (up to an isomorphism)
on the choice of 
$t \in O_{(w^1 + w^2 -r, r)} \;$.

\begin{proposition}\label{GL2SpaltensteinHall}
If the inequality
\begin{equation}\nonumber
r^1 + r^2 \leq r \leq
\min ( w^2 - r^2 + r ^1 , w^1 - r^1 + r^2 )
\end{equation}
holds then the variety 
$\mathfrak{S}_2 (((w^1, r^1), (w^2, r^2)), 
(w^1 + w^2,r))$ is
a smooth connected quasi-projective variety
of dimension
\begin{multline}\nonumber
\dim \mathfrak{S}_2 (((w^1, r^1), (w^2, r^2)), 
(w^1 + w^2,r)) 
= \\ = 
w^1 w^2 + \frac{1}{2} (
\dim O_{(w^1 - r^1 , r^1)}
+ \dim O_{(w^2 - r^2 , r^2)} 
- \dim O_{(w^1 + w^2 - r , r)} ) \; .
\end{multline}
Otherwise 
$\mathfrak{S}_2 (((w^1, r^1), (w^2, r^2)), 
(w^1 + w^2,r))$
is empty.
\end{proposition}
\begin{proof}
Fix $t \in 
\mathfrak{N}_2 (w^1 + w^2)$, $\rank t =r$.
The operator $t$ defines a flag
\begin{equation}\nonumber 
\{ 0 \} \subset \im t \subset \ker t \subset
\mathbb{C}^{w^1 + w^2} \; .
\end{equation}
The variety 
$\mathfrak{S}_2 (((w^1, r^1), (w^2, r^2)), 
(w^1 + w^2,r))$ is a variety of all subspaces 
$X \subset \mathbb{C}^{w^1 + w^2}$
such that 

$\dim (X \cap \im t) = r - r^2$,

$\im ( t |_X ) \subset (X \cap \im t )$,

$\dim ( \im t|_X ) = r^1$,

$\dim ((X \cap \ker t ) / (X \cap \im t )) =
w^1 -r -r^1 +r^2$.

\noindent
The idea of the proof of the Proposition
is to construct a subspace $X$ satisfying the
above conditions in several steps. 
On each step the choices are parametrized 
by a certain variety. Start with
$X \cap \im t$. It is an arbitrary subspace
of dimension $r-r^2$ in $\im t$. So one has 
a Grassmannian $Gr_{r-r^2}^r$ of
such subspaces. Having 
fixed $X \cap \im t$ choose
a subspace $\im ( t |_X )$ inside of it 
(the choices are parametrized by 
the Grassmannian $Gr_{r^1}^{r-r^2}$). 
Next step is to select $X \cap \ker t$. The
intersection $X \cap \im t$ is already
fixed. Thus one proceeds in two steps.
First choose a subspace 
$(X \cap \ker t ) / (X \cap \im t )$
inside $\ker t / \im t $ 
(which gives  
$Gr_{w^1-r-r^1+r^2}^{w^1 + w^2 - 2r}$).
Then each choice of
$X \cap \ker t$ corresponds to an
element of the affine space
$\Hom ((X \cap \ker t ) / (X \cap \im t ) ,
(\im t / (X \cap \im t )))$.
At this point 
$(X \cap \ker t ) \subset \ker t$ and 
$(X / (X \cap \ker t)) \subset 
\mathbb{C}^{w^1 + w^2}/\ker t$ are fixed
(the latter is determined by
$\im ( t |_X )$). Thus the 
remaining choices for $X$
correspond to elements of
$\Hom ((X / (X \cap \ker t )),
(\ker t / (X \cap \ker t )))$.

A rigorous way to spell out
the above argument is to say that one has
the following chain of
locally trivial fibrations
\begin{equation}\nonumber
pt  \leftarrow
A_1 \leftarrow
A_2 \leftarrow
A_3 \leftarrow
A_4 \leftarrow
\mathfrak{S}_2 (((w^1, r^1), (w^2, r^2)), 
(w^1 + w^2,r)) 
\end{equation}
with fibers isomorphic to  
$Gr_{r-r^2}^r$, 
$Gr_{r^1}^{r-r^2}$, 
$Gr_{w^1-r- r^1+r^2}^{w^1 + w^2 - 2r}$, 
$\Hom (\mathbb{C}^{ w^1 - r - r^1 + r^2 } , 
\mathbb{C}^{r^2})$, and
$\Hom (\mathbb{C}^{ r^1 } ,
\mathbb{C}^{w^2 - r + r ^1})$
respectively.

The proposition follows.
\end{proof}

The statement about the dimension of
$\mathfrak{S}_2 (((w^1, r^1), (w^2, r^2)), 
(w^1 + w^2,r))$ is a special case
of the Spaltenstein theorem 
(cf. \ref{Spaltenstein}).

Let
$\mathcal{S}_2 (((w^1, r^1), (w^2, r^2)), 
(w^1 + w^2,r))$ denote
the set of irreducible 
components of the variety 
$\mathfrak{S}_2 (((w^1, r^1), (w^2, r^2)), 
(w^1 + w^2,r))$. Thus
\begin{multline}
\label{GL2S2}
\mathcal{S}_2 (((w^1, r^1), (w^2, r^2)), 
(w^1 + w^2,r))
= \\ =
\begin{cases} 
\{ \mathfrak{S}_2 (((w^1, r^1), (w^2, r^2)), 
(w^1 + w^2,r)) \} 
&\text{if } r \geq r^1 + r^2 \text{ and }\\
&r \leq
\min \{ w^2 \! - \! r^2 \! + \! r^1 , 
w^1 \! - \! r^1 \! + \! r^2 \} 
\\
\emptyset &\text{otherwise } \; .
\end{cases}
\end{multline}
This description of the set of irreducible
components of
$\mathfrak{S}_2 (((w^1, r^1), (w^2, r^2)), 
(w^1 + w^2,r))$ as a one-element or
the empty set is a special case of the
Hall theorem (cf. \ref{HallTheorem}).

Turning back to the variety 
$\mathfrak{T}_2 (w^1, r^1, w^2, r^2)$
recall (\ref{GL2FibredProduct})
that the fiber of the map 
\begin{equation}\nonumber
\delta_2:
\mathfrak{T}_2 (v, w^1, r^1, w^2, r^2)
\rightarrow 
\mathfrak{N}_2 (w^1 + w^2)
\end{equation} 
over a point in $O_{(w^1+w^2-r,r)}$ 
is isomorphic to
\begin{equation}\nonumber
\mathfrak{M}_2 (v, w^1+w^2, r) \times
\mathfrak{S}_2 (((w^1, r^1), (w^2, r^2)), 
(w^1 + w^2,r)) \; .
\end{equation}
Proposition 
\ref{GL2SpaltensteinHall} and
formulas 
\eqref{GL2DimM} and \eqref{GL2DimO2}
imply that 
$\delta_2^{-1} 
(O_{(w^1+w^2-r,r)})$ is a smooth 
connected quasi-projective variety 
of dimension
\begin{multline}\nonumber
\dim \delta_2^{-1} 
(O_{(w^1+w^2-r,r)}) 
= \\ =
w^1 w^2 + v (w^1 + w^2 -v)
+ \frac{1}{2} (
\dim O_{(w^1 - r^1 , r^1)}
+ \dim O_{(w^2 - r^2 , r^2)} ) \; .
\end{multline}
In particular the dimension does not
depend on $r$ (or, in other words, on
the stratum $O_{(w^1+w^2-r,r)}$). 
Therefore the variety 
$\mathfrak{T}_2 (v, w^1, r^1, w^2, r^2)$
is of pure dimension,
\begin{multline}\label{GL2dimT}
\dim \mathfrak{T}_2 (v, w^1, r^1, w^2, r^2) 
= \\ =
w^1 w^2 + v (w^1 + w^2 -v)
+ \frac{1}{2} (
\dim O_{(w^1 - r^1 , r^1)}
+ \dim O_{(w^2 - r^2 , r^2)} ) \; ,
\end{multline}
and the closures of 
$\delta_2^{-1} 
(O_{(w^1+w^2-r,r)})$ are
irreducible components of
the variety 
$\mathfrak{T}_2 (v, w^1, r^1, w^2, r^2)$.
In this way one obtains bijections
\begin{multline}\nonumber
\qquad
\mathcal{T}_2 (v, w^1, r^1, w^2, r^2)
\leftrightarrow
\\
\leftrightarrow
\bigsqcup_{r \in \mathbb{Z}_{\geq 0}}
\mathcal{S}_2 
(((w^1,r^1),(w^2,r^2)),(w^1 + w^2, r))
\times
\mathcal{M}_2 (v, w^1 + w^2, r) \; ,
\end{multline}
\begin{multline}\label{GL2TSMbijection}
\mathcal{T}_2 (w^1, r^1, w^2, r^2)
\leftrightarrow
\\
\leftrightarrow 
\bigsqcup_{r \in \mathbb{Z}_{\geq 0}}
\mathcal{S}_2 
(((w^1,r^1),(w^2,r^2)),(w^1+ w^2, r))
\times
\mathcal{M}_2 (w^1 + w^2, r) \; ,
\end{multline}
where
$\mathcal{T}_2 (v, w^1, r^1, w^2, r^2)$
(resp. $\mathcal{T}_2 (w^1, r^1, w^2, r^2)$)
is the set of irreducible components of
$\mathfrak{T}_2 (v, w^1, r^1, w^2, r^2)$
(resp. $\mathfrak{T}_2 (w^1, r^1, w^2, r^2)$).
These bijections represent the first way
of labeling the set
$\mathcal{T}_2 (w^1, r^1, w^2, r^2)$.
Another labeling of the same set is described 
in the next subsection.

\subsection{The tensor product diagram}
\label{GL2TensorDiagram}
Consider the following 
commutative diagram:
\begin{equation}\label{GL2TTMMGlobal}
\xymatrix{ 
\mathfrak{T}_2 
(v, w^1, r^1, w^2, r^2) 
\ar[d]^{\eta_2}
\ar@/_50pt/[dd]_{\delta_2}   
&
\mathfrak{T}'_2 
(v, w^1, r^1, w^2, r^2) 
\ar[d]^{\eta'_2} 
\ar@{.>}[r]^-{b_2}  
\ar[l]_-{a_2}
&
\mathfrak{M}'_2 (w^1, r^1) 
\! \times \! 
\mathfrak{M}'_2 (w^2, r^2) 
\ar[d]^{\pi_2 \times \pi_2} 
\\
\mathfrak{R}_2 (w^1, r^1, w^2, r^2)
\ar[d]^{\varkappa_2} 
&
\mathfrak{R}'_2 (w^1, r^1, w^2, r^2) 
\ar[r]^-{d_2} 
\ar[l]_-{c_2}  
&
O_{(w^1- r^1, r^1)} \! \times \! 
O_{(w^2- r^2, r^2)} 
\\
\mathfrak{N}_2 ( w^1 + w^2 )
}
\end{equation}
Here the notation is as follows:

$\mathfrak{R}_2 (w^1, r^1, w^2, r^2)$
is a variety of pairs
$(t,X)$, where 
$t \in \mathfrak{N}_2 (w^1+w^2)$,
and $X$ is a subspace of 
$\mathbb{C}^{w^1 + w^2}$
such that 
$\dim X = w^1$, $t X \subset X$,
$\rank t|_X = r^1$,
$\rank t|_{(\mathbb{C}^{w^1+w^2}/X)}
=r^2$;

$\mathfrak{T}'_2 
(v, w^1, r^1, w^2, r^2)$
(resp. $\mathfrak{R}'_2 
(w^1, r^1, w^2, r^2)$)
is the set of
triples $(x , R , Q)$, where 
$x = (t,X,F) \in 
\mathfrak{T}_2 
(v, w^1, r^1, w^2, r^2)$
(resp. $x=(t,X) \in
\mathfrak{R}_2 
(w^1, r^1, w^2, r^2)$), 
$R$ is an isomorphism
$X \rightarrow \mathbb{C}^{w^1}$,
$Q$ is an isomorphism
$\mathbb{C}^{w^1 +w^2}/X 
\rightarrow \mathbb{C}^{w^2}$;

$\mathfrak{M}'_2 (w^i, r^i ) =
\pi_2^{-1} (O_{( w^i-r^i, r^i )})
\subset
\mathfrak{M}_2 ( w^i )$
for $i=1,2$;

$\eta_2 (( t, X, F )) = ( t, X )$;

$\eta'_2 ((x, R, Q)) = 
( \eta_2 (x), R, Q )$;

$\varkappa_2 (( t, X )) = t$;

$a_2 (((t, X, F), R, Q )) = (t, X, F)$;

$c_2 (((t, X), R, Q)) = (t, X)$; 

$b_2 (((t, X, F), R, Q)) \! = \! 
((R t|_{X} R^{-1}, R (F \cap X)),
(Q t|_{\sss (\mathbb{C}^{w^1+w^2} /X)} Q^{-1}, 
Q (F / (F \cap X))))$

$d_2 (((t, X), R, Q)) = 
( R t|_{X} R^{-1} ,
Q t|_{\sss (\mathbb{C}^{w^1+w^2} /X)} Q^{-1})$.

\noindent 
Note that the map $b_2$ is not regular
(this is why the dotted 
arrow is used in the diagram). 
However restricting $a_2$ and $b_2$ to
\begin{equation}\nonumber
a_2^{-1} (\mathfrak{T}_2 
(v^1, w^1, r^1, v^2, w^2, r^2)) =
b_2^{-1}
(\mathfrak{M}'_2 (v^1, w^1, r^1) 
\times 
\mathfrak{M}'_2 (v^2, w^2, r^2))
\end{equation}
one obtains a commutative diagram
in which all maps are regular:
\begin{equation}\label{GL2TTMMRestricted}
\xymatrix{ 
\mathfrak{T}_2 
(v^1, w^1, r^1, v^2, w^2, r^2) 
\ar[d]^{\eta_2}
\ar@/_50pt/[dd]_{\delta_2}   
&
\mathfrak{T}'_2 
(v^1, w^1, r^1, v^2, w^2, r^2) 
\ar[d]^{\eta'_2} 
\ar[r]^-{b_2}  
\ar[l]_-{a_2}
&
*\txt{$\mathfrak{M}'_2 (v^1, w^1, r^1)
\times$\\$
\mathfrak{M}'_2 (v^2, w^2, r^2)\; \;$} 
\ar[d]^{\pi_2 \times \pi_2} 
\\
\mathfrak{R}_2 (w^1, r^1, w^2, r^2)
\ar[d]^{\varkappa_2} 
&
\mathfrak{R}'_2 (w^1, r^1, w^2, r^2) 
\ar[r]^-{d_2} 
\ar[l]_-{c_2}  
&
*\txt{$O_{(w^1- r^1, r^1)} 
\times$\\$ 
O_{(w^2- r^2, r^2)}\; \;$} 
\\
\mathfrak{N}_2 ( w^1 + w^2 )
}
\end{equation}

\noindent Here the notation is as follows

$\mathfrak{T}_2 
(v^1, w^1, r^1, v^2, w^2, r^2)$ is
a locally closed subset of the tensor 
product variety defined in 
subsection \ref{GL2DefinitionOfT};

$\mathfrak{T}'_2 
(v^1, w^1, r^1, v^2, w^2, r^2)=
a_2^{-1} (\mathfrak{T}_2 
(v^1, w^1, r^1, v^2, w^2, r^2))$;

$\mathfrak{M}'_2 (v^i, w^i, r^i ) =
\mathfrak{M}'_2 (w^i, r^i ) \cap 
\mathfrak{M}_2 ( v^i, w^i )$
for $i=1,2$.

\noindent Let
\begin{equation}\nonumber
\begin{split}
\mathfrak{R}_2 (w^1, r^1, w^2, r^2, r) &=
\varkappa_2^{-1} 
(O_{( w^1 + w^2 - r , r )}) \; ,
\\
\mathfrak{R}'_2 (w^1, r^1, w^2, r^2, r) &=
q_1^{-1} 
(\mathfrak{R}_2 
(w^1, r^1, w^2, r^2, r)) \; .
\end{split}
\end{equation}
\begin{proposition}\indent\par
\begin{alphenum}
\item\label{GL2TTMMkappa}
in \eqref{GL2TTMMGlobal} and 
\eqref{GL2TTMMRestricted} 
the map $\varkappa_2$ restricted to
$\mathfrak{R}_2 
(w^1, r^1, w^2, r^2, r)$ is a
locally trivial fibration over
$O_{( w^1 + w^2 - r , r )}$ with a
fiber isomorphic to the Spaltenstein 
variety
$\mathfrak{S}_2 (
((w^1,r^1),(w^2,r^2)),
(w^1+w^2,r))$;
\item\label{GL2TTMMp1q1}
in \eqref{GL2TTMMGlobal}
and \eqref{GL2TTMMRestricted}  
the maps $a_2$ and $c_2$ are principal 
$GL (w^1) \times GL (w^2)$-
bundles;
\item\label{GL2TTMMp2}
in \eqref{GL2TTMMRestricted}
the map $b_2$ is a locally trivial 
fibration with a constant 
smooth connected fiber of dimension
\begin{equation}\nonumber
\dim GL(w^1) +
\dim GL(w^2) +
w^1 w^2 + v^1 (w^2 - v^2)
+v^2 (w^1 - v^1) \; ;
\end{equation}
\item\label{GL2TTMMq2}
if $r^1 + r^2 \leq r \leq
\min ( w^2 - r^2 + r ^1 , w^1 - r^1 + r^2 )$
then  
$d_2$ restricted to 
$\mathfrak{R}'_2 (w^1, r^1, w^2, r^2, r)$
is a locally trivial
fibration with a constant 
smooth connected fiber of dimension
\begin{multline}\nonumber
\qquad \qquad \quad
\dim GL(w^1) +
\dim GL(w^2) 
+ \\ + 
w^1 w^2 +
\frac{1}{2}
(\dim O_{( w^1 + w^2 - r , r)} -
\dim O_{( w^1 - r^1 , r^1)} -
\dim O_{( w^2 - r^2 , r^2)} ) \; ,
\end{multline}
otherwise 
$\mathfrak{R}'_2 (w^1, r^1, w^2, r^2, r)$
is empty;
\item\label{GL2TTMMrank}
if  $r=\min \{ w^2 - v^2 + r^1 , v^1 + r^2 \}$ 
then
$\eta^{\prime -1}_2 (\mathfrak{R}'_2 
(w^1, r^1, w^2, r^2, r))$ 
is an open dense subset of 
$\mathfrak{T}'_2 
(v^1, w^1, r^1, v^2, w^2, r^2)$;
\item\label{GL2TTMMTprime}
$\mathfrak{T}'_2 
(v^1, w^1, r^1, v^2, w^2, r^2)$
is a smooth
connected quasi-projective variety,
\begin{multline}\nonumber
\qquad \qquad \quad
\dim \mathfrak{T}'_2 
(v^1, w^1, r^1, v^2, w^2, r^2)=
\dim GL(w^1) + \dim GL(w^2) 
+ \\ +
 w^1 w^2 + v (w^1 + w^2 -v) +
\frac{1}{2}
(
\dim O_{( w^1 - r^1 , r^1)} +
\dim O_{( w^2 - r^2 , r^2)} ) \; ,
\end{multline}
where $v = v^1 + v^2$;
\item\label{GL2TTMMRprime}
if 
$r^1 + r^2 \leq r \leq
\min ( w^2 - r^2 + r ^1 , w^1 - r^1 + r^2 )$
then  
$\mathfrak{R}'_2 (w^1, r^1, w^2, r^2, r)$
is a smooth
connected quasi-projective variety,
\begin{multline}\nonumber
\qquad \qquad \quad
\dim \mathfrak{R}'_2 
(w^1, r^1, w^2, r^2, r)=
\dim GL(w^1) + \dim GL(w^2) 
+ \\ + 
w^1 w^2 +
\frac{1}{2}
(\dim O_{( w^1 + w^2 - r , r)} +
\dim O_{( w^1 - r^1 , r^1)} +
\dim O_{( w^2 - r^2 , r^2)} ) \; ,
\end{multline}
otherwise 
$\mathfrak{R}'_2 (w^1, r^1, w^2, r^2, r)$
is empty.
\end{alphenum}
\end{proposition}
\begin{proof}
Statements 
\ref{GL2TTMMkappa} and
\ref{GL2TTMMp1q1} 
follow from definitions.
For \ref{GL2TTMMp2} 
write $b_2$ as a composition
of fibrations (cf. the proof of
Proposition \ref{GL2SpaltensteinHall}):
\begin{multline}\nonumber
\mathfrak{T}'_2 
(v^1, w^1, r^1, v^2, w^2, r^2 )  
\xrightarrow{p_1}
\mathfrak{T}''_2 
(v^1, w^1, r^1, v^2, w^2, r^2 )  
\xrightarrow{p_2}
\\
\xrightarrow{p_2}
\mathfrak{T}'''_2
(v^1, w^1, r^1, v^2, w^2, r^2 )
\xrightarrow{p_3}
\mathfrak{T}''''_2
(v^1, w^1, r^1, v^2, w^2, r^2 )
\xrightarrow{p_4}
\\
\xrightarrow{p_4}
\mathfrak{M}'_2
(v^1, w^1, r^1)
\times
\mathfrak{M}'_2 
(v^2, w^2, r^2) \; .
\end{multline}
Here the notation is as follows:

$\mathfrak{T}''''_2 
(v^1, w^1, r^1, v^2, w^2, r^2 )$ 
is a variety of tuples
$(X, t^1, F^1, t^2, F^2)$,
where 
\begin{equation}\nonumber
(( t^1, F^1 ), ( t^2 , F^2))
\in
\mathfrak{M}'_2
(v^1, w^1, r^1) \times
\mathfrak{M}'_2 
(v^2, w^2, r^2)
\end{equation}
and $X$ is a subspace of
$\mathbb{C}^{w^1 + w^2}$
of dimension $w^1$;

$\mathfrak{T}'''_2 
(v^1, w^1, r^1, v^2, w^2, r^2 )$ 
is a variety of tuples
$(R, Q, X, t^1, F^1, t^2, F^2)$,
where 
\begin{equation}\nonumber
(X, t^1, F^1, t^2, F^2)
\in
\mathfrak{T}''''_2 
(v^1, w^1, r^1, v^2, w^2, r^2 )
\end{equation} 
and
$R$ (resp. $Q$) is an isomorphism
$X \rightarrow \mathbb{C}^{w^1}$
(resp . $\mathbb{C}^{w^1+w^2}/X 
\rightarrow \mathbb{C}^{w^2}$);

$\mathfrak{T}''_2 
(v^1, w^1, r^1, v^2, w^2, r^2 )$ 
is a variety of tuples
$(F, R, Q, X, t^1, F^1, t^2, F^2)$,
where 
\begin{equation}\nonumber
(R, Q, X, t^1, F^1, t^2, F^2)
\in
\mathfrak{T}'''_2 
(v^1, w^1, r^1, v^2, w^2, r^2 )
\end{equation} 
and
$F$ is a subspace of 
$\mathbb{C}^{w^1+w^2}$
such that
$R (F \cap X )= F^1$,
$Q (F /(F \cap X)) = F^2$;

\begin{multline}\nonumber
p_1 (((t, X, F), R, Q)) 
= \\ = 
(F, R, Q, X, 
R t|_{X} R^{-1} ,
R (F \cap X) ,
Q t|_{(\mathbb{C}^{w^1+w^2}/X)} Q^{-1},
Q (F/(F \cap X)) ) \; ;
\end{multline}

$p_2$, $p_3$, $p_4$ are natural projections.

The fiber of $p_4$ is a Grassmannian of
dimension $w^1 w^2$, the fiber
of $p_3$ is isomorphic to 
$GL( w^1 )\times GL( w^2 )$,
the fiber of $p_2$ over 
$(R, Q, X, t^1, F^1, t^2, F^2)$
is isomorphic to the affine space 
$\Hom (Q^{-1} ( F^2 ), X / ( R^{-1} (F^1)))$,
and one can describe the fiber of 
$p_1$ over
$(F, R, Q, X, t^1, F^1, t^2, F^2)$
as follows.
Choose a subspace 
$(F \cap X)^X \subset X$
(resp.
$(F \cap X)^F \subset F$,
$(F \cup X)^{\mathbb{C}^{w^1 + w^2}} 
\subset \mathbb{C}^{w^1 + w^2}$)
complimentary to $F \cap X$ in $X$
(resp. 
to $F \cap X$ in $F$,
to $F \cup X$ 
in $\mathbb{C}^{w^1 + w^2}$), 
and consider $R$ (resp. $Q$) as a map 
$(F \cap X)^X \oplus (F \cap X)
\rightarrow \mathbb{C}^{w^1}$
(resp.
$(F \cup X)^{\mathbb{C}^{w^1 + w^2}}
\oplus (F \cap X)^F 
\rightarrow \mathbb{C}^{w^2}$). 
Then
\begin{equation}\nonumber
R^{-1} t^1 R =
\begin{pmatrix}
0 & 0 \\
u^1 & 0 
\end{pmatrix} \; , \qquad
Q^{-1} t^2 Q =
\begin{pmatrix}
0 & 0 \\
u^2 & 0 
\end{pmatrix} \; ,
\end{equation}
for some matrices $u^1$, $u^2$, and
the fiber of $p_1$ consists of operators
$t \in \End ( \mathbb{C}^{w^1 + w^2} )$ 
that have the following block form
\begin{equation}\label{GL2ts}
t=
\begin{pmatrix}
0 & 0 & 0 & 0 \\
0 & 0 & 0 & 0 \\
u^2 & 0 & 0 & 0 \\
s & u^1 & 0 & 0 
\end{pmatrix}
\end{equation}
with respect to the presentation of
$\mathbb{C}^{w^1 + w^2}$ as a direct sum
\begin{equation}\nonumber
\mathbb{C}^{w^1 + w^2} =
(F \cup X)^{\mathbb{C}^{w^1 + w^2}} \oplus
(F \cap X)^X \oplus
(F \cap X)^F \oplus
(F \cap X) \; .
\end{equation}
In \eqref{GL2ts} $s$ is an
arbitrary $v^1 \times (w^2 - v^2)$ matrix.
The statement \ref{GL2TTMMp2} follows. 
Proof of \ref{GL2TTMMq2} is analogues.
The crucial step is to describe
the set of matrices of the form
\begin{equation}\nonumber
t =
\begin{pmatrix}
Q^{-1} t^2 Q & 0 \\
* & R^{-1} t^1 R 
\end{pmatrix} \; ,
\end{equation}
such that $t \in O_{(w^1 + w^2-r,r)}$.

Statement \ref{GL2TTMMrank} follows from
the fact that matrices of rank equal to
$\min (w^2 - v^2 + r^1 , v^1 + r^2 )$ form
an open dense subset in the set of matrices
of the form \eqref{GL2ts}. 

Statement \ref{GL2TTMMTprime}
follows from \ref{GL2TTMMp2} and
\ref{GL2DefinitionOfM}.

Statement \ref{GL2TTMMRprime}
follows from \ref{GL2TTMMq2} or
\ref{GL2TTMMkappa} together with
\ref{GL2TTMMp1q1}.

\end{proof}

Statement \ref{GL2TTMMp1q1} together with 
\eqref{GL2dimT} implies
that the variety 
$\mathfrak{T}'_2 
(v, w^1, r^1, w^2, r^2)$
is of pure dimension,
\begin{multline}\nonumber
\dim \mathfrak{T}'_2 (v, w^1, r^1, w^2, r^2) =
\dim GL(w^1) + \dim GL(w^2) 
+ \\ +
w^1 w^2 + v (w^1 + w^2 -v)
+ \frac{1}{2} (
\dim O_{(w^1 - r^1 , r^1 )} 
+\dim O_{(w^2 - r^2 , r^2 )} ) \; ,
\end{multline}
and the set of irreducible components of
$\mathfrak{T}'_2 
(v, w^1, r^1, w^2, r^2)$ is in a natural
bijection with
$\mathcal{T}_2 
(v, w^1, r^1, w^2, r^2)$ 
(the set of irreducible components of
$\mathfrak{T}_2 
(v, w^1, r^1, w^2, r^2)$.

The fiber of
$\eta'_2$ (in the diagram
\eqref{GL2TTMMGlobal}) over a point in
$\mathfrak{R}'_2
(w^1, r^1, w^2, r^2, r)$ is isomorphic
to $\mathfrak{M}_2 (v, w^1 + w^2 , r)$.
Since 
\begin{equation}\nonumber
\dim \mathfrak{R}'_2 (w^1, r^1, w^2, r^2, r) 
+ \dim \mathfrak{M}_2 (v, w^1 + w^2, r) =
\dim \mathfrak{T}'_2 (v, w^1, r^1, w^2, r^2)
\end{equation}
one obtains another description of the 
bijection \eqref{GL2TSMbijection}
\begin{multline}\nonumber
\mathcal{T}_2 (v, w^1, r^1, w^2, r^2)
\leftrightarrow 
\\
\leftrightarrow
\bigsqcup_{r \in \mathbb{Z}_{\geq 0}}
\mathcal{S}_2 
(((w^1,r^1),(w^2,r^2)),(w^1 + w^2, r))
\times
\mathcal{M}_2 (v, w^1 + w^2, r) \; ,
\end{multline}
where now $\mathcal{S}_2 
(((w^1,r^1),(w^2,r^2)),(w^1 + w^2, r))$
represents the set of irreducible 
components of 
$\mathfrak{R}'_2 (w^1, r^1, w^2, r^2, r)$
(cf. \ref{GL2TTMMRprime})
or the set of irreducible components of a
fiber of $d_2$ in the diagram
\eqref{GL2TTMMGlobal} (cf. \ref{GL2TTMMq2}).
This description the bijection
\eqref{GL2TSMbijection} as coming from the
double fibration $d_2 \circ \eta'_2$ 
is convenient for study of 
the $gl(2)$-restriction in $gl_{\sss N}$
and Kac-Moody cases.

On the other hand 
$\mathfrak{T}'_2 (v, w^1, r^1, w^2, r^2)$
is the total space of the double fibration
$(\pi_2 \times \pi_2 ) \circ b_2$. To make 
sense of $b_2$ one should consider
the diagram \ref{GL2TTMMRestricted} instead
of \ref{GL2TTMMGlobal}. Then 
\ref{GL2TTMMTprime} says that 
$\mathfrak{T}'_2 
(v_1, w^1, r^1, v^2, w^2, r^2)$
is smooth and have the same dimension as
$\mathfrak{T}'_2 (v, w^1, r^1, w^2, r^2)$.
Hence the closure of
$\mathfrak{T}'_2 
(v_1, w^1, r^1, v^2, w^2, r^2)$
is an irreducible component of
$\mathfrak{T}'_2 (v, w^1, r^1, w^2, r^2)$,
and thus one obtains the 
following bijections
\begin{equation}\label{GL2TMMbijection}
\begin{split}
\mathcal{T}_2 
(v^1, w^1, r^1, v^2, w^2, r^2)
&\leftrightarrow
\mathcal{M}_2 
(v^1, w^1, r^1) \times
\mathcal{M}_2 
(v^2, w^2, r^2) \; ,
\\
\mathcal{T}_2 
(v, w^1, r^1, w^2, r^2)
&\leftrightarrow
\bigsqcup_{\substack{
v^1, v^2 \in \mathbb{Z}_{\geq 0}
\\ v^1 + v^2 = v}}
\mathcal{M}_2 
(v^1, w^1, r^1) \times
\mathcal{M}_2 
(v^2, w^2, r^2) \; ,
\\
\mathcal{T}_2 
(w^1, r^1, w^2, r^2)
&\leftrightarrow
\mathcal{M}_2 
(w^1, r^1) \times
\mathcal{M}_2 
(w^2, r^2) \; .
\end{split}
\end{equation}

Combining \eqref{GL2TMMbijection} and
\eqref{GL2TSMbijection} (which represent
two ways of labeling of the set of 
irreducible components of the tensor product 
variety 
$\mathfrak{T}_2 (w^1, r^1, w^2, r^2)$)
one obtains a bijection
\begin{multline}\nonumber
\tau_2 :
\mathcal{M}_2 
(w^1, r^1) \times
\mathcal{M}_2 
(w^2, r^2) 
\xrightarrow{\sim} \\
\xrightarrow{\sim} 
\bigsqcup_{r \in \mathbb{Z}_{\geq 0}}
\mathcal{S}_2 
(((w^1,r^1),(w^2,r^2)),(w^1+ w^2, r))
\times
\mathcal{M}_2 (w^1 + w^2, r) \; .
\end{multline}
It follows from the definition of
$\tau_2$ and \ref{GL2TTMMrank} that 
\begin{multline}\label{GL2TauLabels}
\tau_2 
\bigl( (\mathfrak{M}_2 
(v^1, w^1, r^1) \; , \;
\mathfrak{M}_2 
(v^2, w^2, r^2)) \bigr)
= \\ = 
\bigl( \mathfrak{S}_2 
(((w^1,r^1),(w^2,r^2)),(w^1+ w^2, r_0)) \; , \;
\mathfrak{M}_2 (v^1 + v^2 , w^1 + w^2, r_0) 
\bigr) \; ,
\end{multline}
where
$r_0 = \min \{ w^2 - v^2 + r^1 , v^1 + r^2 \}$.

\subsection{$gl_2$ theorem}

The equality \eqref{GL2TauLabels}
together with definitions 
\eqref{DefinitionOfCrystalProduct} and
\eqref{GL2CrystalDefinition}
imply the following theorem,
which is the main theorem of this paper in
$gl_2$ case.

\begin{theorem}\label{GL2Theorem}
The map $\tau_2$ (described in sections
\ref{GL2FibredProduct}--\ref{GL2TensorDiagram} 
using two labelings 
of the set of irreducible components of
the tensor product variety 
$\mathfrak{T}_2 (w^1, r^1, w^2, r^2)$ is
a crystal isomorphism
\begin{multline}\nonumber
\tau_2 :
\mathcal{M}_2 
(w^1, r^1) \otimes
\mathcal{M}_2 
(w^2, r^2) 
\xrightarrow{\sim} \\
\xrightarrow{\sim} 
\bigoplus_{r \in \mathbb{Z}_{\geq 0}}
\mathcal{S}_2 
(((w^1,r^1),(w^2,r^2)),(w^1+ w^2, r))
\otimes
\mathcal{M}_2 (w^1 + w^2, r) \; ,
\end{multline}
where $\mathcal{S}_2 
(((w^1,r^1),(w^2,r^2)),(w^1+ w^2, r))$ 
is considered as a
set with the trivial $gl_2$-crystal structure.
\end{theorem}

\section{Spaltenstein varieties 
and $gl_{\sss N}$-crystals}
\label{GLnSection}

This section contains a geometric
description of a closed family of
$gl_{\sss N}$-crystals.
Throughout the section $N$ is a fixed positive 
integer.

\subsection{Notation}
Let $\mathcal{Q}_m$ be the set of $m$-tuples of
non-negative integers. Given
$\mathbf{a} \in \mathcal{Q}_m$ let 
$\mathbf{a}_i$ denote the $i$-th component of 
$\mathbf{a}$, and  
$| \mathbf{a} |$ denote the sum of all
components of $\mathbf{a}$:
$| \mathbf{a} | =
\sum_{i=1}^m \mathbf{a}_i$.
Let $\mathcal{Q}_m (k) = 
\{ \mathbf{a} \in \mathcal{Q}_m \; | \;
| \mathbf{a} | = k \}$.

In this paper representations of $gl_{\sss N}$ 
are assumed integrable (polynomial). Therefore all 
weights are positive and one can identify
the relevant part of the weight lattice 
of $gl_{\sss N}$ with $\mathcal{Q}_{\sss N}$. 

Let $\mathcal{Q}_m^+$ (resp. 
$\mathcal{Q}_m^+ (k)$)
be the set of $\mathbf{a} \in \mathcal{Q}_m$
(resp.  $\mathbf{a} \in \mathcal{Q}_m (k)$) 
such that
$\mathbf{a}_1 \geq \mathbf{a}_2 \geq \dots
\geq \mathbf{a}_m$. One can think
of $\mathcal{Q}_m^+$ as of the set of
partitions, or Young diagrams, or
highest weights of integrable
representations of $gl_m$. 
Let $L (\lambda)$ be the irreducible 
highest weight
representation with the highest weight 
$\lambda \in \mathcal{Q}_{\sss N}^+$,
$s_{\lambda}$ be its character 
(a Schur function),
and $\mathcal{L} (\lambda)$ be
the crystal of its canonical basis 
(a highest weight $gl_{\sss N}$-crystal). 

\subsection{Nilpotent orbits
and Spaltenstein varieties}
\label{Spaltenstein}

Let $m \in \mathbb{Z}_{>0}$, 
$w \in \mathbb{Z}_{\geq 0}$,
$t \in \End (\mathbb{C}^w)$,
$t^m = 0$.
Denote by $J(t) \in \mathcal{Q}_m^+ (w)$
a partition given as follows:
$J(t)_i$ is equal to the number of
Jordan blocks of $t$ with size
greater or equal to $i$
($J(t)$ is called below 
the Jordan type of $t$).
For $\lambda \in
\mathcal{Q}_m^+$ let 
\begin{equation}\nonumber
O_{\lambda}=
\{ t \in \End (\mathbb{C}^{| \lambda |}) 
\; | \; t^m = 0, \; J(t) = \lambda \} \; .
\end{equation}
The set $O_{\lambda}
\subset \End (\mathbb{C}^{| \lambda|})$ 
forms a single $GL(| \lambda |)$-orbit.

The stabilizer in $GL(| \lambda |)$ 
of any $t \in O_{\lambda}$ is connected
because it is given by the compliment
of a hypersurface $\det x =0$ in the
affine space $\{ x \in \End 
(\mathbb{C}^{| \lambda |}) \; | \;
xt =tx \}$. Hence $O_{ \lambda }$
is simply connected. It follows that
given a locally trivial fibration over 
$O_{ \lambda }$ with a constant fiber of
pure dimension there is a canonical bijection
between the sets of irreducible components
of the fiber and of the total space of
the fibration. Below these two sets are
usually denoted by the same symbol.

Let $m \in \mathbb{Z}_{>0}$, 
$\mathbf{w}\in \mathcal{Q}_m$. 
Denote by $\mathfrak{F}_m (\mathbf{w})$ 
the partial flag variety
\begin{equation}\nonumber
\mathfrak{F}_m (\mathbf{w})=
\{\mathbf{F}=( \{ 0 \} =
\mathbf{F}_0 \subset \mathbf{F}_1
\subset \ldots \subset 
\mathbf{F}_m=\mathbb{C}^{|\mathbf{w}|} )
\; | \; 
\dim ( \mathbf{F}_{i} / \mathbf{F}_{i-1} ) = 
\mathbf{w}_i \} \; .
\end{equation}
and let $\mathfrak{F}_m (w) = 
\bigsqcup_{\mathbf{w} \in \mathcal{Q}_m (w)}
\mathfrak{F}_m (\mathbf{w})$
be the variety of all
$m$-step partial flags in  
$\mathbb{C}^w$.
Given $\mathbf{F} \in \mathfrak{F}_m (w)$ let
$\dim \mathbf{F}$ denote the $m$-tuple of 
integers consisting of dimensions of subfactors 
of the flag $\mathbf{F}$. Thus 
$\mathbf{F} \in \mathfrak{F}_m (\mathbf{w})$
if and only if  $\dim \mathbf{F}=\mathbf{w}$.

Let $l, m \in \mathbb{Z}_{>0}$, 
$\mu \in \mathcal{Q}^+_m$,
$\Lambda$ be an $l$-tuple 
$(\lambda^1, \dots , \lambda^l)$, where 
$\lambda^i \in 
\mathcal{Q}^+_m $ and
$\sum_{i = 1}^{l}
| \lambda^i | = | \mu |$. Let
$| \Lambda |$ denote
an element of $\mathcal{Q}_l$
given by $| \Lambda | =
( | \lambda^1 | , \dots , | \lambda^l | )$.
Choose $t \in O_{\mu}$ and
consider the following subvariety
of $\mathfrak{F}_l (|\Lambda|)$:
\begin{equation}\nonumber
\mathfrak{S}_l (\Lambda, \mu ) =
\{\mathbf{F} \in \mathfrak{F}_l
(| \Lambda |) 
\; | \;
t \mathbf{F}_i \subset \mathbf{F}_i, \; 
t|_{\mathbf{F}_{i}/\mathbf{F}_{i-1}} 
\in O_{\lambda^i} \} \; .
\end{equation}
The variety $\mathfrak{S}_l (\Lambda, \mu )$
does not depend (up to an isomorphism) on
the choice of $t \in O_{\mu}$. However
when the choice of $t$ is important
the notation is 
$\mathfrak{S}_l (\Lambda, t )$
instead of
$\mathfrak{S}_l (\Lambda, \mu )$.

\begin{theorem}\cite[II.5]{Spaltenstein1982}
The variety 
$\mathfrak{S}_l (\Lambda , \mu )$
is of pure dimension,
\begin{equation}\nonumber
\begin{split}
2 \dim \mathfrak{S}_l (\Lambda, \mu ) 
&= | \mu | (| \mu | - 1)  - 
\sum_{i=1}^l
| \lambda^i | (| \lambda^i | - 1) 
- \dim O_{\mu} + \sum_{i=1}^l
\dim O_{\lambda^i}
= \\ &= 
\sum_{\substack{ i,j = 1 \\ i \neq j}}^l 
| \lambda^i | | \lambda^j |
- \dim O_{\mu} + \sum_{i=1}^l
\dim O_{\lambda^i} \; . 
\end{split}
\end{equation}
\end{theorem}

The variety $\mathfrak{S}_l (\Lambda , \mu )$
is called ($l$-step) Spaltenstein variety.
Let $\mathcal{S}_l (\Lambda , \mu )$
denote the set of irreducible components of
$\mathfrak{S}_l (\Lambda , \mu )$.

\subsection{Hall Theorem}\label{HallTheorem}
Let $\mu^1, \mu^2 \in 
\mathcal{Q}_{\sss N}^+$.
One has the following (tensor product) 
decompositions in the ring of symmetric 
functions, the category of integrable
representations of $gl_{\sss N}$, and
the category of $gl_{\sss N}$-crystals:
\begin{equation}\nonumber
\begin{split}
s_{\mu^1} s_{\mu^2} &=
\sum_{\lambda \in 
\mathcal{Q}_{\sss N}^+ 
(| \mu^1 | + | \mu^2 |)}
c_{\mu^1 \mu^2}^{\lambda} s_{\lambda} \; , 
\\
L (\mu^1) \otimes L (\mu^2) &=
\bigoplus_{\lambda \in 
\mathcal{Q}_{\sss N}^+ 
(| \mu^1 | + | \mu^2 |)}
C ((\mu^1 , \mu^2) , \lambda) 
\otimes L (\lambda)  \; ,
\\
\mathcal{L} (\mu^1) 
\otimes \mathcal{L} (\mu^2) 
&=
\bigoplus_{\lambda \in 
\mathcal{Q}_{\sss N}^+ 
(| \mu^1 | + | \mu^2 |)}
\mathcal{C} ((\mu^1 , \mu^2) , \lambda) 
\otimes \mathcal{L} (\lambda) \; ,
\end{split}
\end{equation}
where $c_{\mu^1 \mu^2}^{\lambda}$ is a
non-negative integer, called the 
Littlewood-Richardson coefficient,
$C ((\mu^1 , \mu^2) , \lambda)$ is
a linear space of dimension
$c_{\mu^1 \mu^2}^{\lambda}$
with the trivial $gl_{\sss N}$-action, and 
$\mathcal{C} ((\mu^1 , \mu^2) , \lambda)$ 
is a trivial $gl_{\sss N}$-crystal (i.e.
a set with the trivial crystal structure) 
of cardinality $c_{\mu^1 \mu^2}^{\lambda}$.

The following theorem is 
due to Hall \cite{Hall1959}.

\begin{theorem}
The number of irreducible components
of the (2-step) Spaltenstein variety
$\mathfrak{S}_2 ((\mu^1,\mu^2), \lambda)$
is equal to the Littlewood-Richardson
coefficient $c_{\mu^1 \mu^2}^{\lambda}$.
In other words, the set
$\mathcal{S}_2 ((\mu^1,\mu^2), \lambda)$ 
is isomorphic to the set
$\mathcal{C} ((\mu^1 , \mu^2) , \lambda)$
as a trivial $gl_{\sss N}$-crystal.
\end{theorem}

\begin{remark}
Let $\mathbb{F}_q$ be a finite field with
$q$ elements. The original Hall's statement  
\cite{Hall1959} says that the number of 
$\mathbb{F}_q$-rational points 
of the variety
$\mathfrak{S}_2 ((\mu^1,\mu^2), \lambda)$
is given by a polynomial in $q$ with 
the leading coefficient equal to 
$c_{\mu^1 \mu^2}^{\lambda}$ 
(the variety 
$\mathfrak{S}_2 ((\mu^1,\mu^2), \lambda)$
is defined over any field).
The above formulation of the theorem
follows from this one by a weak form of
the Weil conjectures 
(cf. \cite{LangWeil1954}). 
\end{remark}

\subsection{A special case}

The following lemma is a special case
of the Hall theorem.

\begin{lemma}\label{GLNSumLemma}
If $\lambda = \mu^1 + \mu^2$ 
(i.e. $\lambda_i = \mu^1_i + \mu^2_i$)
then the variety 
$\mathfrak{S}_2 ((\mu^1,\mu^2), \lambda)$
is non-empty.
\end{lemma}  
\begin{proof}
Let $t \in \mathfrak{N}_{\sss N} 
(| \mu^1 |+| \mu^2 |)$, 
$J (t) = \mu^1 + \mu^2$. 
Choose a basis 
$\{ e_i^j \}_{1 \leq i \leq N}^{
1 \leq j \leq \mu^1_i + \mu^2_i}$
in $\mathbb{C}^{| \mu^1 |+| \mu^2 |}$
such that
\begin{equation}\nonumber
t e_i^j =
\begin{cases}
e_{i-1}^j &\text{ if } i>1 \; , \\
0         &\text{ if } i=1 \; ,
\end{cases}
\end{equation}
and consider a subspace
$X \subset 
\mathbb{C}^{| \mu^1 |+| \mu^2 |}$
defined as the linear span of the set
\begin{multline}\nonumber
\{
e_i^j \quad | \quad 
1 \leq i \leq N \; , 
\\
j \in \mathbb{Z} \cap 
\bigl(
( 0, \mu^1_{\sss N} ] 
\cup 
( \mu^1_{\sss N} + \mu^2_{\sss N} , 
\mu^1_{{\sss N}-1} + \mu^2_{\sss N} ]
\cup \ldots \cup 
( \mu^1_{i+1} + \mu^2_{i+1} , 
\mu^1_{i} + \mu^2_{i+1} ]
\bigr) \} \; ,
\end{multline}
where it is assumed that 
$\mu^1_{{\sss N}+1}=\mu^2_{{\sss N}+1}=0$.
Then the flag 
$(\{ 0 \} \subset X \subset 
\mathbb{C}^{| \mu^1 |+| \mu^2 |})$
belongs to 
$\mathfrak{S}_2 ((\mu^1,\mu^2),
\mu^1 + \mu^2 )$.
\end{proof}

One of the purposes of this paper is
to explain appearance of the varieties 
$\mathfrak{S}_2 ((\mu^1,\mu^2), \lambda)$
in the tensor product decomposition for 
$gl_{\sss N}$. As an application the Hall 
theorem is deduced in \ref{GLNCorollary}
from the Lemma \ref{GLNSumLemma} and 
the Theorem \ref{JosephTheorem}.
The starting point is 
a geometric construction of the highest 
weight crystal $\mathcal{L} (\lambda)$
due to Ginzburg.

\subsection{Ginzburg's construction}
\cite{Ginzburg1991}
Given $w \in \mathbb{Z}_{\geq 0}$ and
$\mathbf{v} \in 
\mathcal{Q}_{\sss N} (w)$
consider the following varieties

$\mathfrak{N}_{\sss N} (w) = 
\{ t \in \End (\mathbb{C}^w) \; | \;
t^N=0 \}$,

$\mathfrak{M}_{\sss N} (\mathbf{v}) 
= \{ (t , \mathbf{F}) \; | \;
t \in \mathfrak{N}_{\sss N} (w) , 
\; \mathbf{F} \in 
\mathfrak{S}_{\sss N} 
( \mathbf{0}, t) \}$, 
where $\mathbf{0}=(0, \dots , 0)$,
$| \mathbf{0} | = \mathbf{v}$
(in other words, 
$\mathfrak{M}_{\sss N} (\mathbf{v})$
is the variety of pairs 
$(t , \mathbf{F})$, where
$t$ is a nilpotent operator 
and $\mathbf{F}$ is an
$N$--step partial flag of
dimension $\mathbf{v}$, such that
$t\mathbf{F}_i \subset \mathbf{F}_{i-1}$),

$\mathfrak{M}_{\sss N} (w)=
\bigsqcup_{\mathbf{v} \in 
\mathcal{Q}_{\sss N} (w)}
\mathfrak{M}_{\sss N} (\mathbf{v})$.

\noindent
Consider a natural map
\begin{equation}\nonumber
\xymatrix{
\mathfrak{M}_{\sss N} (w)
\ar[d]^{\pi_{\sss N}} \\
\; \; \mathfrak{N}_{\sss N} (w) \; , }
\end{equation}
given by $\pi_{\sss N} ((t, \mathbf{F})) = t$,
and put $\mathfrak{M}_{\sss N} (w, t) = 
\pi_{\sss N}^{-1} (t)$. 
One has
\begin{equation}\nonumber
\mathfrak{M}_{\sss N} (w, t) = 
\bigsqcup_{\mathbf{v} \in 
\mathcal{Q}_{\sss N} (w)}
\mathfrak{M}_{\sss N} 
(\mathbf{v}, t) \; ,
\end{equation}
where
\begin{equation}\nonumber
\mathfrak{M}_{\sss N} 
(\mathbf{v}, t) = 
\mathfrak{M}_{\sss N} 
(w, t) \cap 
\mathfrak{M}_{\sss N} 
(\mathbf{v}) \; .
\end{equation}
The variety 
$\mathfrak{M}_{\sss N} (\mathbf{v}, t)$ 
is isomorphic to the Spaltenstein variety
$\mathfrak{S}_{\sss N} (\mathbf{0},t)$.
It depends (up to an isomorphism)
only on $\lambda = J(t)$ 
(the Jordan form of 
$t$, cf. \ref{Spaltenstein}).
Thus the notation
$\mathfrak{M}_{\sss N} (\mathbf{v}, \lambda )$
is often used instead of 
$\mathfrak{M}_{\sss N} (\mathbf{v}, t)$. 
According to \ref{Spaltenstein}
the variety 
$\mathfrak{M}_{\sss N} (\mathbf{v}, 
\lambda )$ is of pure dimension and
\begin{equation}
\label{GLNDimensionOfM}
\dim \mathfrak{M}_{\sss N} 
( \mathbf{v}, \lambda ) =
\frac{1}{2} (\sum_{i \neq j} 
\mathbf{v}_i \mathbf{v}_j  - 
\dim O_{\lambda}) \; .
\end{equation}
Let $\mathcal{M}_{\sss N} 
(\mathbf{v}, \lambda)$ 
be the set of irreducible components of 
$\mathfrak{M}_{\sss N} 
(\mathbf{v}, \lambda )$
and $\mathcal{M}_{\sss N} (\lambda) =
\bigsqcup_{\mathbf{v} \in 
\mathcal{Q}_{\sss N} ( | \lambda | )}
\mathcal{M}_{\sss N} 
(\mathbf{v}, \lambda)$
be the set of irreducible components
of $\mathfrak{M}_{\sss N} 
(| \lambda |, \lambda )$.

The set $\mathcal{M}_{\sss N} (\lambda )$ 
can be endowed with a structure of a normal 
$gl_{\sss N}$-crystal as follows.
The weight function is given by
\begin{equation}\label{GLNWeight}
wt (Z) = \mathbf{v} 
\quad \text{ if } \quad
Z \in \mathcal{M}_{\sss N} 
(\mathbf{v}, \lambda) \; .
\end{equation}
The functions $\varepsilon_k$ and
$\varphi_k$, and the action of the
Kashiwara operators
$\tilde{e}_k$ and $\tilde{f}_k$ 
on $\mathcal{M}_{\sss N} (\lambda )$
are defined in the next subsection using
$gl_2$-restriction.

\subsection{Restriction to $gl_2$}
\label{GLnGL2Section}
Throughout this subsection 
$k$ denotes a fixed integer such that
$1 \leq k \leq N-1$.

Let $\mathcal{Q}_{\sss N}^k$ be the set
of $(N-1)$--tuples of non-negative integers
labeled as follows: $(\mathbf{u}_1,
\dots , \mathbf{u}_{k-1}, \mathbf{u}_{k+1},
\dots , \mathbf{u}_{\sss N})$.
There is a natural map 
$\rho_{\sss N}^k : 
\mathcal{Q}_{\sss N} \rightarrow 
\mathcal{Q}^k_{\sss N}$ given by
\begin{equation}\nonumber
\rho_{\sss N}^k (\mathbf{v}) = 
(\mathbf{v}_1 , \dots , \mathbf{v}_{k-1}, 
\mathbf{v}_k + \mathbf{v}_{k+1},
\dots , \mathbf{v}_{\sss N}) \; .
\end{equation}
Let $\mathcal{Q}^k_{\sss N} (w) =
\rho_{\sss N}^k 
( \mathcal{Q}_{\sss N} (w))=
\{ \mathbf{u} \in \mathcal{Q}_{\sss N}^k
\; | \; \sum \mathbf{u}_i = w \}$. 

Given 
$\mathbf{u} \in \mathcal{Q}_{\sss N}^k$ 
consider the following partial flag variety
\begin{multline}\nonumber
\mathfrak{F}_{\sss N}^k (\mathbf{u})=
\{\mathbf{F}=(
\{ 0 \} =\mathbf{F}_0 \subset \mathbf{F}_1
\subset \ldots \subset 
\mathbf{F}_{k-1} \subset \mathbf{F}_{k+1}
\subset \ldots \subset
\mathbf{F}_{\sss N}
=\mathbb{C}^{|\mathbf{u}|} )
\; | \; \\ 
\dim ( \mathbf{F}_{i} / \mathbf{F}_{i-1} ) = 
\mathbf{u}_i
\text{ for $i \neq 0, k, k+1$ }, \;
\dim ( \mathbf{F}_{k+1} / \mathbf{F}_{k-1})= 
\mathbf{u}_{k+1} \}  \;  ,
\end{multline}
and, similarly, the following variants of
the varieties 
$\mathfrak{M}_{\sss N} (\mathbf{v})$
and $\mathfrak{M}_{\sss N} (w)$:
\begin{multline}\nonumber
\mathfrak{M}_{\sss N}^k (\mathbf{u}) = 
\{ (t , \mathbf{F}) \; | \;
t \in \mathfrak{N}_{\sss N} 
(|\mathbf{u}|) , \; 
\mathbf{F} \in 
\mathfrak{F}_{\sss N}^k (\mathbf{u}) \; , 
\\
t \mathbf{F}_{i} \subset \mathbf{F}_{i-1}
\text{ for $i \neq 0, k, k+1$ }, \;
t \mathbf{F}_{k+1} \subset 
\mathbf{F}_{k+1} \; , \;
t^2 \mathbf{F}_{k+1} \subset 
\mathbf{F}_{k-1}
\} \; ,
\end{multline}
\begin{equation}\nonumber
\mathfrak{M}_{\sss N}^k (w)=
\bigsqcup_{\mathbf{u} \in 
\mathcal{Q}_{\sss N}^k (w)}
\mathfrak{M}_{\sss N}^k (\mathbf{u}) \; .
\end{equation}
One has the following commutative diagram:
\begin{equation}\label{GLnGL2MMN}
\xymatrix{
&\mathfrak{M}_{\sss N} (\mathbf{v}) 
\ar[dl]_{\sigma_{\sss N}^k}
\ar[dd]^{\pi_{\sss N}}
\\
\mathfrak{M}_{\sss N}^k 
(\rho_{\sss N}^k (\mathbf{v}))
\ar[dr]_{\pi_{\sss N}^k}
\\
& \; 
\mathfrak{N}_{\sss N} (|\mathbf{v}|) \; .
}
\end{equation}
Here 
$\sigma_{\sss N}^k$ 
forgets the $k$-th subspace of the flag, and
$\pi_{\sss N}^k$ 
forgets the rest of the flag.

Given 
$\mathbf{u} \in \mathcal{Q}_{\sss N}^k$,
$\mu \in \mathcal{Q}_{\sss N}^+ 
(|\mathbf{u}|)$, and $r \in \mathbb{Z}_{\geq 0}$,
choose $t \in O_{\mu}$ and let
\begin{equation}\nonumber
\begin{split}
\mathfrak{M}_{\sss N}^k (\mathbf{u}, \mu) 
&=
( \pi_{\sss N}^k )^{-1} (t)  \; ,
\\
\mathfrak{M}_{\sss N}^k 
(\mathbf{u}, \mu, r) &= 
\{ (t , \mathbf{F}) 
\in \mathfrak{M}_{\sss N}^k 
(\mathbf{u}, \mu)
\; | \;
\rank (t|_{\mathbf{F}_{k+1} / \mathbf{F}_{k-1}}) 
= r \}  \; .
\end{split}
\end{equation}
The varieties 
$\mathfrak{M}_{\sss N}^k 
(\mathbf{u}, \mu )$ and
$\mathfrak{M}_{\sss N}^k 
(\mathbf{u}, \mu, r)$
do not depend (up to an isomorphism)
on $t \in O_{\mu}$. 
The latter one is an example of $(N-1)$-step
Spaltenstein variety. More exactly, 
$\mathfrak{M}_{\sss N}^k 
(\mathbf{u}, \mu, r) =
\mathfrak{S}_{{\sss N}-1} (\Lambda, \mu )$,
where $\Lambda = 
(0, \dots, 0, \lambda^k , 0 , \dots 0)$ with
$\lambda^k=(\mathbf{u}_{k+1}-r , r)$.
Hence $\mathfrak{M}_{\sss N}^k 
(\mathbf{u},\mu , r)$ is of pure dimension
and
\begin{equation}\label{GLnGL2DimOfM}
\dim \mathfrak{M}_{\sss N}^k 
(\mathbf{u},\mu , r) = 
\frac{1}{2} (
\sum^N_{\substack{
i,j = 1 \\ 
i,j \neq k\\
i \neq j}}
\mathbf{u}_i  \mathbf{u}_j  
- \dim O_{\mu} + 
\dim O_{(\mathbf{u}_{k+1} - r, r)}) \; . 
\end{equation} 
Let $\mathcal{M}_{\sss N}^k 
(\mathbf{u}, \mu , r)$ 
be the set of irreducible components of
$\mathfrak{M}_{\sss N}^k 
(\mathbf{u},\mu , r)$.

One has the following stratification
\begin{equation}\nonumber
\mathfrak{M}_{\sss N}^k 
(\mathbf{u}, \mu)=
\bigcup_{r \in \mathbb{Z}_{\geq 0}}
\mathfrak{M}_{\sss N}^k 
(\mathbf{u}, \mu ,r) 
\end{equation}
with only finitely many non-empty strata.
The map $\sigma_{\sss N}^k$ 
(cf. \eqref{GLnGL2MMN}) restricted to
a stratum
$\mathfrak{M}_{\sss N}^k 
(\rho_{\sss N}^k (\mathbf{v}),\mu ,r)$
is a locally trivial fibration with
the fiber isomorphic to the 
fiber of the map  
\begin{equation}\nonumber
\xymatrix{
\mathfrak{M}_2 (\mathbf{v}_k +
\mathbf{v}_{k+1})
\ar[d]^{\pi_{2}}
\\
\mathfrak{N}_2 (\mathbf{v}_k + 
\mathbf{v}_{k+1}) 
}
\end{equation}
over a point in 
$O_{(\mathbf{v}_{k+1} + 
\mathbf{v}_{k}-r , r)}$
(see \eqref{GL2DefinitionOfM} for
the relevant definitions). 
In other words the fiber is a 
Grassmannian and 
\begin{equation}\nonumber
\dim ((\sigma_{\sss N}^k)^{-1} (x))=
\mathbf{v}_{k} \mathbf{v}_{k+1}  -
\frac{1}{2}(\dim 
O_{(\mathbf{v}_{k} + 
\mathbf{v}_{k+1} -r,r)}) \; ,
\end{equation} 
where $x \in \mathfrak{M}_{\sss N}^k 
(\rho_{\sss N}^k (\mathbf{v}),\mu ,r)$.
Hence the dimension of 
$\mathfrak{M}_{\sss N}^k 
(\rho_{\sss N}^k (\mathbf{v}),\mu ,r)$
plus the dimension of the fiber of
$\sigma_{\sss N}^k$ 
over a point in this stratum 
is equal to the dimension
of $\mathfrak{M}_{\sss N} 
(\mathbf{v},\mu)$ 
(cf. \eqref{GLNDimensionOfM})
and one has
a bijection of the sets of
irreducible components
\begin{equation}\nonumber
\theta_{\sss N}^k :
\bigsqcup_{\mathbf{v}= 
(\rho^k_{\sss N})^{-1} (\mathbf{u})}
\mathcal{M}_{\sss N} (\mathbf{v}, \mu)
\xrightarrow{\sim \;}
\bigsqcup_{r}
\mathcal{M}_{\sss N}^k 
(\mathbf{u}, \mu , r)
\times
\mathcal{M}_2 
(\mathbf{u}_{k+1},r) \; ,
\end{equation}
or
\begin{equation}\label{GLNRestrictionBijection}
\theta_{\sss N}^k :
\mathcal{M}_{\sss N} (\mu)
\xrightarrow{\sim \;}
\bigsqcup_{\substack{\mathbf{u} \in 
\mathcal{Q}_{\sss N}^k 
\\ r \in \mathbb{Z}_{\geq 0}}}
\mathcal{M}_{\sss N}^k 
(\mathbf{u}, \mu , r)
\times
\mathcal{M}_2 
(\mathbf{u}_{k+1},r) \; .
\end{equation}

\subsection{The crystal structure}

Now one can finish the definition
of the crystal structure on the set
$\mathcal{M}_{\sss N} (\mu )$. 
Namely, consider 
the RHS of \ref{GLNRestrictionBijection}
as a $gl_2$-crystal with the crystal structure
coming from the second multiple, and let
\begin{equation}\nonumber
\begin{split}
\varepsilon_k &=  
\varepsilon \circ \theta_{\sss N}^k \; ,
\\
\varphi_k &=  
\varphi \circ \theta_{\sss N}^k \; ,
\\
\tilde{e}_k &= 
(\theta_{\sss N}^k)^{-1} 
(\tilde{e}) (\theta_{\sss N}^k) \; ,
\\
\tilde{f}_k &= 
(\theta_{\sss N}^k)^{-1} 
(\tilde{f}) (\theta_{\sss N}^k) \; .
\end{split}
\end{equation}
These formulas 
together with the weight function
\eqref{GLNWeight} provide a structure 
of $gl_{\sss N}$-crystal on 
$\mathcal{M}_{\sss N} (\mu )$.
By abuse of notation this crystal
is also denoted by 
$\mathcal{M}_{\sss N} (\mu )$.

\begin{proposition}\label{GLNHighestWeight}
The crystal $\mathcal{M}_{\sss N} (\mu )$
is a highest weight normal 
$gl_{\sss N}$-crystal with highest
weight $\mu$.
\end{proposition}
\begin{proof}
The fact that $\mathcal{M}_{\sss N} (\mu )$
is a normal crystal follows immediately
from definitions. To prove that it is 
highest weight it is enough to show
that it contains unique element $Z$ such that
\begin{equation}\label{GLNHWProof1}
\Tilde{e}_k Z = 0 
\quad \text{ for $1 \leq k \leq N-1$ }.
\end{equation}
Fix $t \in \mathfrak{N}_{\sss N} (| \mu |)$ 
such that $J (t) = \mu$, and let 
$(t, \mathbf{F})$ be a generic
point of a component $Z$ of 
$\mathfrak{M}_{\sss N} (t )$ satisfying 
\eqref{GLNHWProof1} (``generic'' here means
``not lying in any other component''). 
By definition of
the operators $\Tilde{e}_k$ it implies that
\begin{equation}\nonumber
\ker t |_{(\mathbf{F}_{k+1}/\mathbf{F}_{k-1})} = 
\mathbf{F}_{k}/\mathbf{F}_{k-1}
\quad \text{ for $1 \leq k \leq N-1$ }.
\end{equation}
From this it follows by induction on $k$ that 
\begin{equation}\nonumber
\ker t^{k-1} |_{\mathbf{F}_{k}} = 
\mathbf{F}_{k-1}
\quad \text{ for $2 \leq k \leq N$ },
\end{equation}
and then using inverse induction on $k$ one 
obtains 
\begin{equation}\nonumber
\mathbf{F}_{k} = \ker t^k 
\quad \text{ for $1 \leq k \leq N$ }.
\end{equation}
Hence the component $Z$ has one generic point
$(t, (\{ 0 \} \subset \ker t \subset \ldots 
\subset \ker t^N = \mathbb{C}^{| \mu |}))$.
It means that there is unique highest weight
component $Z$ of 
$\mathfrak{M}_{\sss N} (| \mu |)$, 
$Z$ is the point, and $wt (Z)= \mu$. 
The proposition follows.
\end{proof}

Actually the crystal 
$\mathcal{M}_{\sss N} ( \mu )$
is isomorphic to
$\mathcal{L} ( \mu )$
(the crystal of the canonical basis
of the irreducible highest weight
module with highest weight $\mu$).
One can prove it using a slightly modified
argument of \cite{Ginzburg1991}. Another
proof (based on the Theorem 
\ref{JosephTheorem}) is given in subsection
\ref{GLNCorollary}.

The above description of the crystal
structure via $gl_2$-restriction 
is essentially due to Lusztig 
\cite[12]{Lusztig1991a}
(in the more general case 
of quiver varieties), 
and was used later by 
Kashiwara and Saito \cite{KashiwaraSaito},
and Nakajima \cite{Nakajima1998}.

\subsection{The tensor product variety and
the tensor product diagram}
\label{GLNTensorSection}
Let $\mu^1 , \mu^2 \in
\mathcal{Q}_{\sss N}^+$.
Consider the following variety:
\begin{multline}\nonumber
\mathfrak{T}_{\sss N} (\mu^1 , \mu^2 ) 
= \\ =
\{ (t, X, \mathbf{F}) \; | \;
t \in \mathfrak{N}_{\sss N} 
(| \mu^1 | + | \mu^2 |), \;
(t,\mathbf{F}) \in \mathfrak{M}_{\sss N} 
(| \mu^1 | + | \mu^2 |), \;
X \in \mathfrak{S}_2 
((\mu^1, \mu^2), t)\} \; .
\end{multline}
In plain words a point of 
$\mathfrak{T}_{\sss N} (\mu^1 , \mu^2 )$ 
is a triple consisting of a nilpotent
operator $t \in \End 
(\mathbb{C}^{|\mu^1|+|\mu^2|})$,
an $N$-step partial flag $\mathbf{F}$ 
in $\mathbb{C}^{|\mu^1|+|\mu^2|}$, 
and a subspace 
$X \subset \mathbb{C}^{|\mu^1|+|\mu^2|}$ 
of dimension $|\mu^1|$, such that
$t$ preserves both $\mathbf{F}$ and
$X$, and when restricted to $X$ 
(resp. $\mathbb{C}^{|\mu^1|+|\mu^2|}/X$,
subfactors of $\mathbf{F}$) gives
a nilpotent operator of Jordan type
$\mu^1$ (resp. 
nilpotent operator of Jordan type $\mu^2$, 
$0$ operator).

The variety
$\mathfrak{T}_{\sss N} (\mu^1 , \mu^2 )$
is called 
\emph{the tensor product variety}
(for $gl_{\sss N}$). It plays
an important role in geometric description 
of the tensor product. 
In particular using two different
ways of labeling irreducible components
of $\mathfrak{T}_{\sss N} (\mu^1 , \mu^2 )$
one can describe the direct sum
decomposition of the tensor product
of $gl_{\sss N}$-crystals
$\mathcal{M}_{\sss N} (\mu^1)$ and
$\mathcal{M}_{\sss N} (\mu^2)$. 
Details are given below.

The variety 
$\mathfrak{T}_{\sss N} (\mu^1 , \mu^2 )$
is a disjoint union of subvarieties
labeled by the dimensions of the
subspaces of $\mathbf{F}$:
\begin{equation}\nonumber
\mathfrak{T}_{\sss N} 
(\mu^1 , \mu^2 )=
\bigsqcup_{\mathbf{v} \in 
\mathcal{Q}_{\sss N} 
(| \mu^1 | + | \mu^2 |)}
\mathfrak{T}_{\sss N} 
(\mathbf{v}, \mu^1 , \mu^2 ) \; ,
\end{equation}
where
\begin{equation}\nonumber
\mathfrak{T}_{\sss N} 
( \mathbf{v}, \mu^1 , \mu^2 )=
\{ (t, X, \mathbf{F}) \in 
\mathfrak{T}_{\sss N} 
( \mu^1 , \mu^2 ) 
\; | \;
\dim \mathbf{F} =
\mathbf{v} \} \; .
\end{equation}
There is a natural map:
\begin{equation}\nonumber
\xymatrix{
\mathfrak{T}_{\sss N} 
( \mathbf{v}, \mu^1 , \mu^2 )
\ar[d]^{\delta_{\sss N}}
\\ \; \; 
\mathfrak{N}_{\sss N} 
(| \mu^1 | + | \mu^2 |) \;,
}
\end{equation}
given by $\delta_{\sss N}
(t, X, \mathbf{F}) = t$. Stratify 
$\mathfrak{N}_{\sss N} 
(| \mu^1 |+| \mu^2 |)$
by $GL(| \mu^1 | + | \mu^2 |)$-orbits.
The map $\delta_{\sss N}$
restricted to an orbit
$O_\lambda$ is a fibration with a constant
fiber isomorphic to the direct
product of two Spaltenstein
varieties (for $2$-step and
$N$-step flags):
\begin{equation}\nonumber
\delta_{\sss N}^{-1} (x) \approx
\mathfrak{S}_2 
(( \mu^1 , \mu^2 ), \lambda ),
\times
\mathfrak{M}_{\sss N} 
( \mathbf{v} , \lambda )
\end{equation}
if $x \in O_{\lambda}$. Therefore
\begin{equation}\nonumber
\dim 
(\delta_{\sss N}^{-1} (O_\lambda))
= | \mu^1 | | \mu^2 |
+ \frac{1}{2}(
\sum_{i \neq j} 
\mathbf{v}_i \mathbf{v}_j
+ \dim O_{\mu^1} + \dim O_{\mu^2}) \; .
\end{equation}
Since the dimension of 
$\delta_{\sss N}^{-1} (O_\lambda)$
does not depend on $\lambda$, it
follows that the variety
$\mathfrak{T}_{\sss N} 
(\mathbf{v}, \mu^1 , \mu^2 )$
is of pure dimension,
\begin{equation}\label{GLnDimensionOfT}
\dim 
\mathfrak{T}_{\sss N} 
(\mathbf{v}, \mu^1 , \mu^2 )
= | \mu^1 | | \mu^2 |
+ \frac{1}{2}(
\sum_{i \neq j} 
\mathbf{v}_i \mathbf{v}_j
+ \dim O_{\mu^1} + \dim O_{\mu^2}) \; ,
\end{equation}
and one has the following bijection
\begin{equation}\label{GLnTSM}
\mathcal{T}_{\sss N} 
(\mathbf{v}, \mu^1 , \mu^2)
\leftrightarrow
\bigsqcup_{\lambda \in 
\mathcal{Q}_{\sss N}^+ 
(| \mu^1 |+| \mu^2 |)}
\mathcal{S}_2 ((\mu^1 , \mu^2), \lambda)
\times
\mathcal{M}_{\sss N} 
(\mathbf{v}, \lambda) \; , 
\end{equation}
where 
$\mathcal{T}_{\sss N} 
(\mathbf{v}, \mu^1, \mu^2 )$ 
is the set of irreducible components of
$\mathfrak{T}_{\sss N} 
(\mathbf{v}, \mu^1 , \mu^2)$.

Another way to label 
irreducible components of 
$\mathfrak{T}_{\sss N} 
(\mu^1 , \mu^2)$
is as follows. Let $\mathbf{F} \cap X$ 
(resp. $\mathbf{F} / (\mathbf{F} \cap X)$) 
denote
the flag in $X$ (resp. in 
$\mathbb{C}^{| \mu^1 |+| \mu^2 |} / X$)
obtained by taking intersections
with $X$ (resp. quotients by intersections
with $X$) of the subspaces of $\mathbf{F}$,
and let
\begin{multline}\nonumber
\mathfrak{T}_{\sss N} 
( \mathbf{v}^1, \mathbf{v}^2, 
\mu^1 , \mu^2 ) 
= \\ =
\{ (t, X, \mathbf{F}) \in 
\mathfrak{T}_{\sss N} 
( \mu^1 , \mu^2 ) 
\; | \;
\dim (\mathbf{F} \cap X) = \mathbf{v}^1, \;
\dim ( \mathbf{F} / (\mathbf{F} \cap X)) = 
\mathbf{v}^2 \} \; . 
\end{multline}
Each subset $\mathfrak{T}_{\sss N} 
( \mathbf{v}^1, 
\mathbf{v}^2, \mu^1 , \mu^2 ) \subset
\mathfrak{T}_{\sss N} 
( \mu^1 , \mu^2 )$
is locally closed and one has
\begin{equation}\nonumber
\mathfrak{T}_{\sss N} 
(\mathbf{v} , \mu^1 , \mu^2 )=
\bigsqcup_{\substack{
\mathbf{v}^1 \in 
\mathcal{Q}_{\sss N} 
(| \mu^1 |) \\
\mathbf{v}^2 \in 
\mathcal{Q}_{\sss N} 
(| \mu^2 |) \\
\mathbf{v}^1 + \mathbf{v}^2 = \mathbf{v}
}}
\mathfrak{T}_{\sss N} 
(\mathbf{v}^1, \mathbf{v}^2 , 
\mu^1 , \mu^2 ) \; .
\end{equation}
Consider the following commutative diagram
(cf. \eqref{GL2TTMMGlobal}):

\begin{equation}\label{GLnTTMMGlobal}
\xymatrix{
\mathfrak{T}_{\sss N} 
(\mathbf{v}, \mu^1 , \mu^2 ) 
\ar[d]^{\eta_{\sss N}}
\ar@/_50pt/[dd]_{\delta_{\sss N}}  
&
\mathfrak{T}'_{\sss N} 
(\mathbf{v}, \mu^1 , \mu^2 )  
\ar[l]_-{a_{\sss N}} 
\ar@{.>}[r]^-{b_{\sss N}}
\ar[d]^{\eta'_{\sss N}}
&
\mathfrak{M}'_{\sss N} ( \mu^1 )
\times
\mathfrak{M}'_{\sss N} ( \mu^2 )
\ar[d]^{\pi_{\sss N} \times \pi_{\sss N}}
\\ 
\mathfrak{R}_{\sss N} ( \mu^1, \mu^2 )
\ar[d]^{\varkappa_{\sss N}}
&
\mathfrak{R}'_{\sss N} ( \mu^1, \mu^2 )
\ar[l]_-{c_{\sss N}} 
\ar[r]^-{d_{\sss N}}
&
\; \; 
O_{\mu^1} \times O_{\mu^2} 
\\
\mathfrak{N}_{\sss N} 
(| \mu^1 |+| \mu^2 |)
}
\end{equation}
Here the notation is as follows:

$\mathfrak{R}_{\sss N} ( \mu^1, \mu^2 )$
is the variety of pairs
$(t,X)$, where 
$t \in \mathfrak{N}_{\sss N} 
(| \mu^1 |+| \mu^2 |)$,
and $X$ is a subspace of 
$\mathbb{C}^{| \mu^1 |+| \mu^2 |}$
such that 
$\dim X = | \mu^1 |$, $t X \subset X$,
$t|_X$ has Jordan type $\mu^1$,
$t|_{(\mathbb{C}^{w^1+w^2}/X)}$
has Jordan type $\mu^2$;

$\mathfrak{T}'_{\sss N} 
(\mathbf{v}, \mu^1 , \mu^2 )$ 
(resp. 
$\mathfrak{R}'_{\sss N} ( \mu^1, \mu^2 )$)
is the variety of triples 
$(x , R , Q)$, where 
$x = (t, X, \mathbf{F}) \in
\mathfrak{T}_{\sss N} 
(\mathbf{v}, \mu^1 , \mu^2 )$
(resp.
$x = (t, X) \in
\mathfrak{R}_{\sss N} ( \mu^1, \mu^2 )$), 
$R$ (resp. $Q$) is an isomorphism
$X \rightarrow \mathbb{C}^{|\mu^1|}$
(resp . $\mathbb{C}^{|\mu^1|+|\mu^2|}/X 
\rightarrow \mathbb{C}^{|\mu^2|}$);

$\mathfrak{M}'_{\sss N} ( \mu^i )=
\pi_{\sss N}^{-1} (O_{\mu^i})
\subset 
\mathfrak{M}_{\sss N} (|\mu^i|)$
for $i=1,2$;

$\eta_{\sss N} (( t, X, \mathbf{F} )) = 
( t, X )$;

$\eta'_{\sss N} ((x, R, Q)) = 
( \eta_{\sss N} (x), R, Q )$;

$\varkappa_{\sss N} (( t, X )) = t$;

$a_{\sss N} 
(((t, X, \mathbf{F}), R, Q )) = 
(t, X, \mathbf{F})$;

$c_{\sss N} (((t, X), R, Q)) = (t, X)$; 

$b_{\sss N} (((t, X, \mathbf{F}), R, Q)) = 
\bigl(\bigl(R t|_{X} R^{-1},
R (\mathbf{F} \cap X)\bigr),
\bigl( 
Q t|_{(\mathbb{C}^{|\mu^1|+|\mu^2|}/X)} Q^{-1}, 
Q (\mathbf{F}/(\mathbf{F} \cap X))
\bigr)\bigr)$;

$d_{\sss N} (((t, X), R, Q)) = 
( R t|_{X} R^{-1} ,
Q t|_{(\mathbb{C}^{w^1+w^2} /X)} Q^{-1})$.

\noindent
The map $b_{\sss N}$ is not regular. 
To fix it one
has to restrict $a_{\sss N}$ and 
$b_{\sss N}$ to the locally closed subset
\begin{equation}\nonumber
\mathfrak{T}'_{\sss N} 
(\mathbf{v}^1, \mathbf{v}^2, \mu^1 , \mu^2 )=
a_{\sss N}^{-1} (\mathfrak{T}_{\sss N} 
(\mathbf{v}^1, \mathbf{v}^2, \mu^1 , \mu^2 ))=
b_{\sss N}^{-1} (
\mathfrak{M}'_{\sss N} ( \mathbf{v}^1, \mu^1 )
\times
\mathfrak{M}'_{\sss N} ( \mathbf{v}^2, \mu^2 )
) \; .
\end{equation}
The result is the following diagram 
(cf. \eqref{GL2TTMMRestricted}:

\begin{equation}\label{GLnTTMMRestricted}
\xymatrix{
\mathfrak{T}_{\sss N} 
(\mathbf{v}^1, \mathbf{v}^2, \mu^1 , \mu^2 ) 
\ar[d]^{\eta_{\sss N}}
\ar@/_50pt/[dd]_{\delta_{\sss N}}  
&
\mathfrak{T}'_{\sss N} 
(\mathbf{v}^1, \mathbf{v}^2, \mu^1 , \mu^2 )  
\ar[l]_-{a_{\sss N}} 
\ar[r]^-{b_{\sss N}}
\ar[d]^{\eta'_{\sss N}}
&
\mathfrak{M}'_{\sss N} 
( \mathbf{v}^1, \mu^1 )
\times
\mathfrak{M}'_{\sss N} 
( \mathbf{v}^2 , \mu^2 )
\ar[d]^{\pi_{\sss N} \times \pi_{\sss N}}
\\ 
\mathfrak{R}_{\sss N} ( \mu^1, \mu^2 )
\ar[d]^{\varkappa_{\sss N}}
&
\mathfrak{R}'_{\sss N} ( \mu^1, \mu^2 )
\ar[l]_-{c_{\sss N}} 
\ar[r]^-{d_{\sss N}}
&
\; \; 
O_{\mu^1} \times O_{\mu^2} 
\\
\mathfrak{N}_{\sss N} 
(| \mu^1 |+| \mu^2 |)
}
\end{equation}

The following proposition is an analogue
of the Proposition \ref{GL2TensorDiagram}.

\begin{proposition}\indent\par
\begin{alphenum}
\item
\label{GLnTTMMp1}
in the diagrams \eqref{GLnTTMMGlobal} and
\eqref{GLnTTMMRestricted} the maps
$a_{\sss N}$ and $c_{\sss N}$ are
principal 
$GL (| \mu^1 |) \times GL (| \mu^2 |)$-bundles;
\item\label{GLnTTMMp2}
in the diagram \eqref{GLnTTMMRestricted}
the map $b_{\sss N}$ 
is a locally trivial fibration 
with a constant smooth connected fiber 
of dimension
\begin{equation}\nonumber
\dim (GL( | \mu^1 | )) +
\dim (GL( | \mu^2 | )) +
| \mu_1 | | \mu_2 | +
\sum_{i\neq j}
\mathbf{v}^1_i \mathbf{v}^2_j \; .
\end{equation}
\end{alphenum}
\end{proposition}
\begin{proof}
Statement \ref{GLnTTMMp1} 
follows from definitions.
To prove \ref{GLnTTMMp2}
write $b_2$ as a composition
of fibrations (cf. the proof of
\ref{GL2TTMMp2}):
\begin{multline}\nonumber
\mathfrak{T}'_{\sss N} 
(\mathbf{v}^1, \mathbf{v}^2, 
\mu^1 , \mu^2 )  
\xrightarrow{q_1}
\mathfrak{T}''_{\sss N} 
(\mathbf{v}^1, \mathbf{v}^2, 
\mu^1 , \mu^2 )  
\xrightarrow{q_2}
\mathfrak{T}'''_{\sss N} 
(\mathbf{v}^1, \mathbf{v}^2, 
\mu^1 , \mu^2 )  
\xrightarrow{q_3}
\\
\xrightarrow{q_3}
\mathfrak{T}''''_{\sss N} 
(\mathbf{v}^1, \mathbf{v}^2, 
\mu^1 , \mu^2 )  
\xrightarrow{q_4}
\mathfrak{M}'_{\sss N} 
(\mathbf{v}^1, \mu^1 )
\times
\mathfrak{M}'_{\sss N} 
(\mathbf{v}^2, \mu^2 ) \; .
\end{multline}
Here the notation is as follows:

$\mathfrak{T}''''_{\sss N} 
(\mathbf{v}^1, \mathbf{v}^2, 
\mu^1 , \mu^2 )$ is a variety of
tuples
$(X,\mathbf{F}^1, \mathbf{F}^2, t^1, t^2)$,
where 
\begin{equation}\nonumber
(( \mathbf{F}^1, t^1 ), ( \mathbf{F}^2, t^2))
\in
\mathfrak{M}'_{\sss N} 
(\mathbf{v}^1, \mu^1 ) \times
\mathfrak{M}'_{\sss N} 
(\mathbf{v}^2, \mu^2 )
\end{equation}
and $X$ is a subspace of
$\mathbb{C}^{| \mu^1 | + | \mu^2 |}$
of dimension $| \mu^1 |$;

$\mathfrak{T}'''_{\sss N} 
(\mathbf{v}^1, \mathbf{v}^2, 
\mu^1 , \mu^2 )$ is a variety of
tuples
$(R, Q, X,
\mathbf{F}^1, \mathbf{F}^2, t^1, t^2)$,
where 
\begin{equation}\nonumber
(X, \mathbf{F}^1, \mathbf{F}^2, t^1, t^2)
\in
\mathfrak{T}''''_{\sss N} 
(\mathbf{v}^1, \mathbf{v}^2, \mu^1 , \mu^2 )
\end{equation} 
and
$R$ (resp. $Q$) is an isomorphism
$X \rightarrow \mathbb{C}^{|\mu^1|}$
(resp . $\mathbb{C}^{|\mu^1|+|\mu^2|}/X 
\rightarrow \mathbb{C}^{|\mu^2|}$);

$\mathfrak{T}''_{\sss N} 
(\mathbf{v}^1, \mathbf{v}^2, 
\mu^1 , \mu^2 )$ is a variety of
tuples
$(\mathbf{F}, R, Q, X, 
\mathbf{F}^1, \mathbf{F}^2, t^1, t^2)$,
where 
\begin{equation}\nonumber
(R, Q, X, 
\mathbf{F}^1, \mathbf{F}^2, t^1, t^2)
\in
\mathfrak{T}'''_{\sss N} 
(\mathbf{v}^1, \mathbf{v}^2, \mu^1 , \mu^2 )
\end{equation} 
and
$\mathbf{F}$ is an $N$-step flag in
$\mathbb{C}^{ | \mu^1 | + | \mu^2 | }$
such that
$R (\mathbf{F} \cap X )= \mathbf{F}^1$,
$Q (\mathbf{F} /(\mathbf{F} \cap X)) = 
\mathbf{F}^2$;

\begin{multline}\nonumber
q_1 (((t, X, \mathbf{F}), R, Q)) = \\ = 
(\mathbf{F},
R,Q,X,
R (\mathbf{F} \cap X) ,
Q (\mathbf{F}/(\mathbf{F} \cap X)) ,
R t|_{X} R^{-1},
Q 
t|_{(\mathbb{C}^{|\mu^1|+|\mu^2|}/X)} 
Q^{-1}) \; ,
\end{multline}

$q_2$, $q_3$, $q_4$ are natural projections.

\noindent
The fiber of $q_4$ is a Grassmannian of
dimension $| \mu^1 ||\mu^2 |$, the fiber
of $q_3$ is isomorphic to 
$GL( | \mu^1 | )\times GL( | \mu_2 | )$,
the fiber of $q_2$ is isomorphic to the affine
space

\begin{equation}\nonumber
\bigoplus_{i} \Hom (
\mathbf{F}^2_i / \mathbf{F}^2_{i-1} , 
\mathbb{C}^{| \mu^1 |} / \mathbf{F}^1_i ) \; ,
\end{equation}
and the fiber of $q_1$ is isomorphic to the 
affine space

\begin{equation}\nonumber
\bigoplus_{i} \Hom (
\mathbf{F}^2_i / \mathbf{F}^2_{i-1} , 
\mathbf{F}^1_{i-1} ) \; .
\end{equation}
Statement \ref{GLnTTMMp2} follows. 

\end{proof}
Statement \ref{GLnTTMMp1} implies 
that the set of irreducible components
of $\mathfrak{T}'_{\sss N} 
(\mathbf{v}, \mu^1 , \mu^2 )$
is in a natural bijection with 
$\mathcal{T}_{\sss N} 
(\mathbf{v}, \mu^1 , \mu^2 )$
( the set of irreducible
components of 
$\mathfrak{T}'_{\sss N} 
(\mathbf{v}, \mu^1 , \mu^2 )$). 
On the other hand \ref{GLnTTMMp2} 
provides a bijection between
the product of the sets of irreducible 
components of 
$\mathfrak{M}_{\sss N} 
(\mathbf{v}^1, \mu^1)$
and 
$\mathfrak{M}_{\sss N} 
(\mathbf{v}^2, \mu^2)$,
and the set of irreducible components of
$\mathfrak{T}'_{\sss N} 
(\mathbf{v}^1, \mathbf{v}^2, 
\mu^1 , \mu^2 )$. 
Moreover it follows from
\ref{GLnTTMMp1}, \ref{GLnTTMMp2},
\ref{GLNDimensionOfM}, and 
\ref{GLnDimensionOfT} that
\begin{equation}\nonumber
\dim \mathfrak{T}'_{\sss N} 
(\mathbf{v}^1, \mathbf{v}^2, 
\mu^1 , \mu^2 ) = 
\dim \mathfrak{T}'_{\sss N} 
(\mathbf{v}^1 + \mathbf{v}^2 , 
\mu^1 , \mu^2 ) \; .
\end{equation}
Hence the closures of
irreducible components of 
$\mathfrak{T}'_{\sss N} 
(\mathbf{v}^1, \mathbf{v}^2, 
\mu^1 , \mu^2 )$ 
are irreducible components of 
$\mathfrak{T}'_{\sss N} 
(\mathbf{v}^1 + \mathbf{v}^2, 
\mu^1 , \mu^2 )$ and one obtains a
bijection
\begin{equation}\nonumber
\mathcal{T}_{\sss N} 
(\mathbf{v} , \mu^1 , \mu^2 )
\leftrightarrow
\bigsqcup_{\substack{
\mathbf{v}^1 \in 
\mathcal{Q}_{\sss N} (| \mu^1 |) \\
\mathbf{v}^2 \in 
\mathcal{Q}_{\sss N} (| \mu^2 |) \\
\mathbf{v}^1 + \mathbf{v}^2 = \mathbf{v}
}}
\mathcal{M}_{\sss N} 
(\mathbf{v}^1, \mu^1 )
\times
\mathcal{M}_{\sss N} 
(\mathbf{v}^2, \mu^2 ) \; .
\end{equation}
Together with \eqref{GLnTSM}
it provides a bijection
\begin{multline}\nonumber
\tau_{\sss N} : \quad
\bigsqcup_{\substack{
\mathbf{v}^1 \in 
\mathcal{Q}_{\sss N} (| \mu^1 |) \\
\mathbf{v}^2 \in 
\mathcal{Q}_{\sss N} (| \mu^2 |) \\
\mathbf{v}^1 + \mathbf{v}^2 = \mathbf{v}
}}
\mathcal{M}_{\sss N} 
(\mathbf{v}^1, \mu^1 )
\times
\mathcal{M}_{\sss N} 
(\mathbf{v}^2, \mu^2 ) 
\xrightarrow{\sim}
\\
\xrightarrow{\sim}
\bigsqcup_{\lambda \in 
\mathcal{Q}_{\sss N}^+ 
(| \mu^1 |+| \mu^2 |)}
\mathcal{S}_2 ((\mu^1 , \mu^2), \lambda)
\times
\mathcal{M}_{\sss N} 
(\mathbf{v}, \lambda ) \; ,
\end{multline}
or, taking union over 
$\mathbf{v} \in \mathcal{Q}_{\sss N} 
(|\mu^1|+|\mu^2|)$,
\begin{equation}\nonumber
\tau_{\sss N} : \quad
\mathcal{M}_{\sss N} (\mu^1 )
\times
\mathcal{M}_{\sss N} (\mu^2 ) 
\xrightarrow{\sim \;}
\bigsqcup_{\lambda \in 
\mathcal{Q}_{\sss N}^+
(| \mu^1 |+| \mu^2 |)}
\mathcal{S}_2((\mu^1 , \mu^2), \lambda)
\times
\mathcal{M}_{\sss N} (\lambda ) \;  .
\end{equation}

\subsection{The main theorem 
for $gl_{\sss N}$}

Here is the main theorem of this paper
in the $gl_{\sss N}$ case.

\begin{theorem}\label{GLnTheorem}
The map $\tau_{\sss N}$ is a crystal 
isomorphism
\begin{equation}\nonumber
\tau_{\sss N} : \quad
\mathcal{M}_{\sss N} (\mu^1 )
\otimes
\mathcal{M}_{\sss N} (\mu^2 ) 
\xrightarrow{\sim}
\bigoplus_{\lambda \in 
\mathcal{Q}_{\sss N}^+ (| \mu^1 |+| \mu^2 |)}
\mathcal{S}_2((\mu^1 , \mu^2), \lambda)
\otimes
\mathcal{M}_{\sss N} (\lambda ) \; ,
\end{equation}
where 
$\mathcal{S}_2((\mu^1 , \mu^2), \lambda)$
is considered as a set 
with the trivial crystal structure.
\end{theorem}

The proof of this theorem occupies the next 
subsection.

\subsection{$gl_2$-restriction 
in the tensor product or the proof
of Theorem \ref{GLnTheorem}}
By definition the bijection 
$\tau_{\sss N}$
preserves the weight function. So in 
order to prove that $\tau_{\sss N}$ is a
crystal morphism
\begin{equation}\nonumber
\tau_{\sss N} :
\mathcal{M}_{\sss N} (\mu^1 )
\otimes
\mathcal{M}_{\sss N} (\mu^2 ) 
\rightarrow
\bigsqcup_{\lambda \in 
\mathcal{Q}_{\sss N}^+ (| \mu^1 |+| \mu^2 |)}
\mathcal{S}_2((\mu^1 , \mu^2), \lambda)
\otimes
\mathcal{M}_{\sss N} (\lambda )  
\end{equation}
it remains to compare the action of
the Kashiwara operators
$\Tilde{f}_k$ and $\Tilde{e}_k$ on the
domain and the range of $\tau_{\sss N}$
(recall that both sides are normal crystals,
and thus one does not need to care about
the functions
$\varepsilon_k$ and $\varphi_k$). 
The idea of the proof is to reduce
the problem to $gl_2$ case
and then use the corresponding result for
$gl_2$ (Theorem \ref{GL2Theorem}).

Let $\mu^1 , \mu^2 \in
\mathcal{Q}_{\sss N}^+$, 
and $k$ be an integer such that
$1 \leq k \leq N-1$.
Consider the following variety
(cf. \ref{GLnGL2Section}):
\begin{multline}\nonumber
\mathfrak{T}^k_{\sss N} (\mu^1 , \mu^2 ) 
= \\ = 
\{ (t, X, \mathbf{F}) \; | \;
t \in \mathfrak{N}_{\sss N} 
(| \mu^1 | + | \mu^2 |), \;
(t,\mathbf{F}) \in \mathfrak{M}^k_{\sss N} 
(| \mu^1 | + | \mu^2 |), \;
X \in \mathfrak{S}_2 
((\mu^1, \mu^2), t)\} \; .
\end{multline}
See \ref{GLnGL2Section} for the
definition of $\mathfrak{M}^k_{\sss N}$.
The variety 
$\mathfrak{T}^k_{\sss N} (\mu^1 , \mu^2 )$
represents $gl_2$ restriction
for the tensor product variety 
$\mathfrak{T}_{\sss N} (\mu^1, \mu^2)$.
One has the following commutative diagram
(cf. \eqref{GLnGL2MMN}):
\begin{equation}\nonumber
\xymatrix{
& \mathfrak{T}_{\sss N} ( \mu^1, \mu^2 ) 
\ar[ld]_{\xi^k_{\sss N}} 
\ar[dd]^{\delta_{\sss N}} 
\\
\mathfrak{T}_{\sss N}^k ( \mu^1, \mu^2 )
\ar[dr]_{\delta^k_{\sss N}} 
\\
& \; \; 
\mathfrak{N}_{\sss N} (|\mu^1|+|\mu^2|)
\; ,}
\end{equation}
where $\xi^k_{\sss N}$
forgets the $k$-th subspace in the 
flag $\mathbf{F}$ and
$\delta^k_{\sss N}$
forgets the rest of the flag. Given 
$\mathbf{u} \in \mathcal{Q}_{\sss N}^k
(| \mu^1 |+| \mu^2 |)$, 
and $r \in \mathbb{Z}_{\geq 0}$
consider the following strata in
$\mathfrak{T}_{\sss N}^k ( \mu^1, \mu^2 )$:
\begin{equation}\nonumber
\mathfrak{T}_{\sss N}^k 
( \mathbf{u}, \mu^1, \mu^2, r ) =
\{ (t, X, \mathbf{F}) \in
\mathfrak{T}_{\sss N}^k ( \mu^1, \mu^2 )
\; | \;
\dim \mathbf{F} = \mathbf{u}, \;
\rank( t|_{\mathbf{F}_{k+1}/\mathbf{F}_{k-1}} ) 
= r \} \; .
\end{equation}
The fiber of 
$\xi^k_{\sss N}$ over a point in 
$\mathfrak{T}_{\sss N}^k 
(\mathbf{u}, \mu^1, \mu^2, r)$
is isomorphic to 
$\mathfrak{M}_2 (\mathbf{u}_{k+1},r)$
(cf. \ref{GL2DefinitionOfM}). 
Consider stratification of 
$\mathfrak{N}_{\sss N} (|\mu^1|+|\mu^2|)$ by
$GL(|\mu^1|+|\mu^2|)$-orbits 
$\{ O_{\lambda} \}$, 
The map $\delta_{\sss N}^k$ restricted to 
$O_{\lambda}$ is a locally trivial
fibration with a fiber isomorphic
to the product of two Spaltenstein varieties
$\mathfrak{S}_2 (( \mu^1 , \mu^2 ), \lambda )$ 
and $\mathfrak{M}_{\sss N}^k 
(\mathbf{u}, \lambda, r)$ 
(cf. \ref{GLnGL2Section}). 
Hence
$\mathfrak{T}_{\sss N}^k 
(\mathbf{u}, \mu^1, \mu^2, r)$ is of pure
dimension,
\begin{multline}\label{GLnGL2DimOfTNk}
\dim \mathfrak{T}_{\sss N}^k 
(\mathbf{u}, \mu^1, \mu^2, r) =
| \mu^1 ||\mu^2 | + \\ +
\frac{1}{2}(
\sum_{i \neq j}
\mathbf{u}_i \mathbf{u}_j
+ \dim O_{\mu^1} + \dim O_{\mu^2}
+ \dim O_{(\mathbf{u}_{k+1}-r , r)}) \; , 
\end{multline}
and one has the following bijection
\begin{equation}\label{GLnGL2TSM}
\mathcal{T}_{\sss N}^k 
(\mathbf{u}, \mu^1, \mu^2, r)
\leftrightarrow
\bigsqcup_{\lambda \in 
\mathcal{Q}_{\sss N}^+ (| \mu^1 |+| \mu^2 |)}
\mathcal{S}_2 ((\mu^1, \mu^2), \lambda)
\times
\mathcal{M}_{\sss N}^k 
(\mathbf{u}, \lambda, r) \; , 
\end{equation}
where 
$\mathcal{T}_{\sss N}^k 
(\mathbf{u}, \mu^1, \mu^2, r)$
is the set of irreducible components of
$\mathfrak{T}_{\sss N}^k 
(\mathbf{u}, \mu^1, \mu^2, r)$.

Recall that the fiber of 
$\xi^k_{\sss N}$ over a point in 
$\mathfrak{T}_{\sss N}^k 
(\mathbf{u}, \mu^1, \mu^2, r)$
is isomorphic to 
$\mathfrak{M}_2 (\mathbf{u}_{k+1},r)$.
Since 
\begin{equation}\nonumber
\dim \mathfrak{T}_{\sss N}^k 
(\rho_{\sss N}^k (\mathbf{v}), 
\mu^1, \mu^2, r) +
\dim \mathfrak{M}_2 
(\mathbf{v}_{k},
\mathbf{v}_k+\mathbf{v}_{k+1}, r ) =
\dim \mathfrak{T}_{\sss N}
( \mathbf{v}, \mu^1, \mu^2 )
\end{equation}
one has a bijection:
\begin{equation}\label{GLnGL2TTM}
\mathcal{T}_{\sss N} (\mu^1, \mu^2)
\leftrightarrow
\bigsqcup_{\substack{
\mathbf{u} \in \mathcal{Q}_{\sss N}^k 
(| \mu^1 |+| \mu^2 |) \\
r \in \mathbb{Z}_{\geq 0}}}
\mathcal{T}_{\sss N}^k 
(\mathbf{u}, \mu^1, \mu^2, r)
\times 
\mathcal{M}_2 
(\mathbf{u}_{k+1}, r) \; .
\end{equation}
Combining \eqref{GLnGL2TSM} with
\eqref{GLnGL2TTM} one obtains a
bijection
\begin{equation}\label{GLnGL2TSMM}
\mathcal{T}_{\sss N} (\mu^1, \mu^2)
\leftrightarrow
\bigsqcup_{\substack{
\lambda \in \mathcal{Q}_{\sss N}^+ 
(| \mu^1 |+| \mu^2 |)\\
\mathbf{u} \in \mathcal{Q}_{\sss N}^k
(| \mu^1 |+| \mu^2 |) \\
r \in \mathbb{Z}_{\geq 0}}}
\mathcal{S}_2 ((\mu^1, \mu^2), \lambda)
\times
\mathcal{M}_{\sss N}^k (\mathbf{u}, \lambda , r)
\times 
\mathcal{M}_2 (\mathbf{u}_{k+1}, r) \; .
\end{equation}
This bijection is the result of 
$gl_2$-restriction
on the range of $\tau_{\sss N}$. In particular,
the operators $\Tilde{e}_k$ and $\Tilde{f}_k$ 
act on the last multiple in the RHS of
\eqref{GLnGL2TSMM}.

To carry out $gl_2$-restriction on 
the domain of $\tau_{\sss N}$ one has to 
forget the $k$-th subspace of
the flag $\mathbf{F}$
in all varieties in the 
diagrams \ref{GLnTTMMGlobal}
and \ref{GLnTTMMRestricted}. 
The results are the following
commutative diagrams
\begin{equation}\label{GLnGL2TTMMGlobal}
\xymatrix{
\mathfrak{T}_{\sss N} (\mathbf{v}, 
\mu^1 , \mu^2)
\ar[dd]_-{\xi^k_{\sss N}} 
&
\mathfrak{T}'_{\sss N} (\mathbf{v}, 
\mu^1 , \mu^2)  
\ar[l]_-{a_{\sss N}^k} 
\ar@{.>}[r]^-{b_{\sss N}^k}
\ar[dd]_-{\xi^{\prime k}_{\sss N}} 
&
\mathfrak{M}'_{\sss N} (\mu^1 )
\times
\mathfrak{M}'_{\sss N} (\mu^2 )
\ar[dd]_-{\sigma^k_{\sss N} 
\times \sigma^k_{\sss N}} 
\\ \\
\mathfrak{T}_{\sss N}^k 
(\rho_{\sss N}^k (\mathbf{v}) , 
\mu^1 , \mu^2 ) 
&
\mathfrak{T}^{\prime k}_{\sss N} 
(\rho_{\sss N}^k (\mathbf{v}) ,
\mu^1 , \mu^2 )  
\ar[l]_-{c_{\sss N}^k} 
\ar@{.>}[r]^-{d_{\sss N}^k}
&
\mathfrak{M}^{\prime k}_{\sss N} 
( \mu^1 )
\times
\mathfrak{M}^{\prime k}_{\sss N}
( \mu^2 )
\ar[dd]_-{\pi_{\sss N}^k \times \pi_{\sss N}^k}
\\ \\ & &
\; \; 
O_{\mu^1} \times O_{\mu^2}
\; ,}
\end{equation}
\noindent and
\begin{equation}\label{GLnGL2TTMMRestricted}
\xymatrix{
\mathfrak{T}_{\sss N} (\mathbf{v}^1, \mathbf{v}^2, 
\mu^1 , \mu^2)
\ar[dd]_-{\xi^k_{\sss N}} 
&
\mathfrak{T}'_{\sss N} (\mathbf{v}^1, \mathbf{v}^2, 
\mu^1 , \mu^2)  
\ar[l]_-{a_{\sss N}^k} 
\ar[r]^-{b_{\sss N}^k}
\ar[dd]_-{\xi^{\prime k}_{\sss N}} 
&
\mathfrak{M}'_{\sss N} (\mathbf{v}^1, \mu^1 )
\times
\mathfrak{M}'_{\sss N} (\mathbf{v}^2, \mu^2 )
\ar[dd]_-{\sigma^k_{\sss N} 
\times \sigma^k_{\sss N}} 
\\ \\
\mathfrak{T}_{\sss N}^k 
(\mathbf{u}^1, \mathbf{u}^2, 
\mu^1 , \mu^2 ) 
&
\mathfrak{T}^{\prime k}_{\sss N} 
(\mathbf{u}^1, \mathbf{u}^2, 
\mu^1 , \mu^2 )  
\ar[l]_-{c_{\sss N}^k} 
\ar[r]^-{d_{\sss N}^k}
&
\mathfrak{M}^{\prime k}_{\sss N} 
(\mathbf{u}^1, \mu^1 )
\times
\mathfrak{M}^{\prime k}_{\sss N}
(\mathbf{u}^2, \mu^2 )
\ar[dd]_-{\pi_{\sss N}^k \times \pi_{\sss N}^k}
\\ \\ & &
\; \; 
O_{\mu^1} \times O_{\mu^2}
}
\end{equation}
where 
$\mathbf{u}^1 = \rho_{\sss N}^k 
( \mathbf{v}^1 )$,
$\mathbf{u}^2 = \rho_{\sss N}^k 
( \mathbf{v}^2 )$.
In the above two diagrams
the top rows are copies of 
the diagrams \ref{GLnTTMMGlobal} and
\ref{GLnTTMMRestricted};

$\mathfrak{T}^{k}_{\sss N} 
(\mathbf{u}^1 , \mathbf{u}^2 , 
\mu^1 , \mu^2 ) =
\{ (t , X ,\mathbf{F}) \in
\mathfrak{T}^{k}_{\sss N} 
( \mu^1 , \mu^2 ) 
\; | \;
\dim (\mathbf{F} \cap X)=\mathbf{u}^1 , \;
\dim(\mathbf{F}/(\mathbf{F}\cap X)) =
\mathbf{u}^2  \}$;

$\mathfrak{T}^{\prime k}_{\sss N} 
(\mathbf{u}^1 , \mathbf{u}^2 , 
\mu^1 , \mu^2 )$ 
(resp. 
$\mathfrak{T}^{\prime k}_{\sss N} 
(\mathbf{u}, \mu^1 , \mu^2 )$)
is the variety of triples
$(x,R,Q)$, where $x=(t, X, \mathbf{F}) \in
\mathfrak{T}^{k}_{\sss N} 
(\mathbf{u}^1 , \mathbf{u}^2 , 
\mu^1 , \mu^2 )$
(resp. $x=(t, X, \mathbf{F}) \in
\mathfrak{T}^{k}_{\sss N} 
(\mathbf{u}, \mu^1 , \mu^2 )$), and
$R$ (resp. $Q$) is an isomorphism
$X \rightarrow \mathbb{C}^{| \mu^1 |}$
(resp. 
$\mathbb{C}^{| \mu^1 |+| \mu^2 |}/X
\rightarrow \mathbb{C}^{| \mu^2 |}$);

the vertical maps are natural projections
(forgetting the $k$-th subspace of
the flag $\mathbf{F}$);

the maps $c_{\sss N}^k$ and 
$d_{\sss N}^k$ are uniquely
defined to make the diagrams commutative.

\noindent
The maps $c_{\sss N}^k$ and 
$d_{\sss N}^k$ are regular only
in the diagram \eqref{GLnGL2TTMMRestricted}.

Let 
\begin{multline}\nonumber
\mathfrak{T}^k_{\sss N} 
(\mathbf{u}^1, \mathbf{u}^2, 
\mu^1 , \mu^2 , r^1, r^2, r) 
= \\ =
\{ (t, X, \mathbf{F}) \in
\mathfrak{T}^k_{\sss N} 
( \mathbf{u}^1, \mathbf{u}^2, \mu^1 , \mu^2 ) 
\; | \; 
\rank (t|_{
\mathbf{F}_{k+1}/ \mathbf{F}_{k-1}}) 
=r , 
\\
\rank (t|_{
(X \cap \mathbf{F}_{k+1})/
(X \cap \mathbf{F}_{k-1})}) 
=r^1 , \;
\rank (t|_{
(\mathbf{F}_{k+1}/((\mathbf{F}_{k+1} \cap X)+
\mathbf{F}_{k-1})}) 
=r^2 \} \; ,
\end{multline}
and 
\begin{equation}\nonumber
\mathfrak{T}^{\prime k}_{\sss N} 
(\mathbf{u}^1, \mathbf{u}^2, 
\mu^1 , \mu^2 , r^1 , r^2 , r) =
(c_{\sss N}^k)^{-1} (
\mathfrak{T}^k_{\sss N} 
(\mathbf{u}^1, \mathbf{u}^2, 
\mu^1 , \mu^2 , r^1, r^2, r) ) \; .
\end{equation}

The set
$\mathfrak{T}^k_{\sss N} 
(\mathbf{u}^1, \mathbf{u}^2, 
\mu^1 , \mu^2 , r^1, r^2, r)$
(resp. 
$\mathfrak{T}^{\prime k}_{\sss N} 
(\mathbf{u}^1, \mathbf{u}^2, 
\mu^1 , \mu^2 , r^1 , r^2 , r)$)
is a locally closed subset in
$\mathfrak{T}^k_{\sss N} 
(\mathbf{u}^1 +  \mathbf{u}^2, 
\mu^1 , \mu^2 , r)$
(resp. 
$\mathfrak{T}^{\prime k}_{\sss N} 
(\mathbf{u}^1 + \mathbf{u}^2, 
\mu^1 , \mu^2 , r)$).

The proof of the following proposition
is similar to the proofs of Propositions
\ref{GL2TensorDiagram}
and \ref{GLNTensorSection}.

\begin{proposition}\indent\par
\begin{alphenum}
\item\label{GLnGL2TTMMp1q1}
in the diagrams \eqref{GLnGL2TTMMGlobal}
and \eqref{GLnGL2TTMMRestricted} the maps
$a_{\sss N}^k$ and $c_{\sss N}^k$ 
are principal
$GL (| \mu^1 |) \times GL (| \mu^2 |)$-
bundles;
\item\label{GLnGL2TTMMq2}
if $r^1 + r^2 \leq r \leq
\min ( \mathbf{u}^2_{k+1} - r^2 + r ^1 , 
\mathbf{u}^1_{k+1} - r^1 + r^2 )$
then $d_{\sss N}^k$ restricted to 
$\mathfrak{T}^{\prime k}_{\sss N} 
(\mathbf{u}^1, \mathbf{u}^2, 
\mu^1 , \mu^2 , r^1 , r^2 , r)$
is a locally trivial
fibration over
$\mathfrak{M}^{\prime k}_{\sss N} 
( \mathbf{u}^1 , \mu^1, r^1 )
\times
\mathfrak{M}^{\prime k}_{\sss N}
( \mathbf{u}^2 , \mu^2, r^2 )$
with a constant 
smooth connected fiber of dimension
\begin{multline}\nonumber
\qquad \qquad \quad
\dim GL(| \mu^1 |) +
\dim GL(| \mu^2 |) 
+ | \mu^1 || \mu^2 |
+\sum_{i \neq j}\mathbf{u}^1_i \mathbf{u}^2_j
+ \\ +
\frac{1}{2}
(\dim O_{( \mathbf{u}^1_{k+1} 
+ \mathbf{u}^2_{k+1} - r , r)} 
-\dim O_{( \mathbf{u}^1_{k+1} - r^1 , r^1)} 
-\dim O_{( \mathbf{u}^2_{k+1} - r^2 , r^2)} ) \; ,
\end{multline}
otherwise 
$\mathfrak{T}^{\prime k}_{\sss N} 
(\mathbf{u}^1, \mathbf{u}^2, 
\mu^1 , \mu^2 , r^1 , r^2 , r)$
is empty.
\end{alphenum}
\end{proposition}

It follows from \ref{GLnGL2TTMMp1q1} 
and \ref{GLnGL2DimOfTNk} that
$\mathfrak{T}^{\prime k}_{\sss N} 
(\mathbf{u}, \mu^1, \mu^2, r)$ is of pure
dimension,
\begin{multline}\label{GLnGL2DimOfTprime}
\dim \mathfrak{T}^{\prime k}_{\sss N} 
(\mathbf{u}, \mu^1, \mu^2, r) =
\dim GL (| \mu^1 |) +
\dim GL (| \mu^2 |) +
| \mu^1 ||\mu^2 | 
+ \\ +
\frac{1}{2} (
\sum_{i \neq j}
\mathbf{u}_i \mathbf{u}_j
+ \dim O_{\mu^1} + \dim O_{\mu^2}
+ \dim O_{(\mathbf{u}_{k+1}-r , r)}) \; , 
\end{multline}
and the set of irreducible components of
$\mathfrak{T}^{\prime k}_{\sss N} 
(\mathbf{u}, \mu^1, \mu^2, r)$ is in a natural
bijection with 
$\mathcal{T}^{k}_{\sss N} 
(\mathbf{u}, \mu^1, \mu^2, r)$
(the set of irreducible components of
$\mathfrak{T}^{k}_{\sss N} 
(\mathbf{u}, \mu^1, \mu^2, r)$).

On the other hand it follows from
\ref{GLnGL2TTMMq2},
\eqref{GLnGL2DimOfM}, and
\eqref{GLnGL2DimOfTprime} 
that the variety 
$\mathfrak{T}^{\prime k}_{\sss N} 
(\mathbf{u}^1, \mathbf{u}^2, 
\mu^1, \mu^2, r^1, r^2, r)$ has the same
dimension as 
$\mathfrak{T}^{\prime k}_{\sss N} 
(\mathbf{u}^1 + \mathbf{u}^2, 
\mu^1, \mu^2, r)$
(in particular, the dimension does not
depend on $r^1$, $r^2$, and 
$\mathbf{u}^1 - \mathbf{u}^2$). Hence
irreducible components of 
$\mathfrak{T}^{\prime k}_{\sss N} 
(\mathbf{u}^1 + \mathbf{u}^2, 
\mu^1, \mu^2, r)$ are closures of
irreducible components of
$\mathfrak{T}^{\prime k}_{\sss N} 
(\mathbf{u}^1, \mathbf{u}^2, 
\mu^1, \mu^2, r^1, r^2, r)$ and using 
\ref{GLnGL2TTMMq2} one obtains the
following bijection:
\begin{multline}\label{GLnGL2TSMMsmall}
\mathcal{T}^k_{\sss N} (
\mathbf{u} , 
\mu^1 , \mu^2 , r) 
\leftrightarrow
\bigsqcup_{\substack{
\mathbf{u}^1 \in 
\mathcal{Q}_{\sss N}^k (| \mu^1 |) \\
\mathbf{u}^2 \in 
\mathcal{Q}_{\sss N}^k (| \mu^2 |) \\
\mathbf{u}^1 + \mathbf{u}^2 = \mathbf{u}\\
r^1, \; r^2 \in \mathbb{Z}_{\geq 0}
}}
\mathcal{S}_2 
(((\mathbf{u}^1_{k+1},r^1), 
(\mathbf{u}^2_{k+1},r^2)),
(\mathbf{u}_{k+1},r)) 
\times \\ \times
\mathcal{M}_{\sss N}^k
(\mathbf{u}^1,\mu^1 , r^1) 
\times
\mathcal{M}_{\sss N}^k
(\mathbf{u}^2, \mu^2 , r^2)
\end{multline} 
where 
$\mathcal{S}_2 
(((\mathbf{u}^1_{k+1},r^1), 
(\mathbf{u}^2_{k+1},r^2)),
(\mathbf{u}_{k+1},r))$ 
comes from the fiber of $d_{\sss N}^k$
(cf. the inequality in \ref{GLnGL2TTMMq2}).

Combining \eqref{GLnGL2TSMMsmall} 
with \eqref{GLnGL2TSM}
one obtains a bijection
\begin{multline}\label{GLnDefinitionOfTauNk}
\tau_{\sss N}^k : \quad 
\bigsqcup_{\substack{
\mathbf{u}^1 \in 
\mathcal{Q}_{\sss N}^k (| \mu^1 |) \\
\mathbf{u}^2 \in 
\mathcal{Q}_{\sss N}^k (| \mu^2 |) \\
\mathbf{u}^1 + \mathbf{u}^2 = \mathbf{u}\\
r^1, \; r^2 \in \mathbb{Z}_{\geq 0}
}}
\mathcal{S}_2 
(((\mathbf{u}^1_{k+1},r^1), 
(\mathbf{u}^2_{k+1},r^2)),
(\mathbf{u}_{k+1},r)) 
\times \\ \times 
\mathcal{M}_{\sss N}^k
(\mathbf{u}^1, \mu^1 , r^1) 
\times 
\mathcal{M}_{\sss N}^k
(\mathbf{u}^2, \mu^2 , r^2)
\xrightarrow{\sim}
\\ \\
\xrightarrow{\sim}
\bigsqcup_{\lambda \in 
\mathcal{Q}_{\sss N}^+ 
(| \mu^1 |+| \mu^2 |)}
\mathcal{S}_2 ((\mu^1, \mu^2), \lambda)
\times
\mathcal{M}_{\sss N}^k 
(\mathbf{u}, \lambda, r) \; .
\end{multline}

Let $\mu^1$, 
$\mu^2 \in \mathcal{Q}^+_{\sss N}$, and
$\mathbf{v} \in \mathcal{Q}_{\sss N}
(| \mu^1 |+| \mu^2 |)$. Consider a diagram

\begin{equation}\label{GLnGL2main}
\vcenter{
\xymatrix{
*\txt{$\bigsqcup$ \\
$\substack{\lambda \in 
\mathcal{Q}^+_{\sss N} 
(| \mu^1 |+| \mu^2 |)}
$}
&
*\txt{$
\mathcal{S}_2 ((\mu^1 , \mu^2 ), \lambda)
\times
\mathcal{M}_{\sss N} (\mathbf{v}, \lambda)
$}
& \\
*\txt{$\bigsqcup$ \\
$\substack{
\mathbf{v}^1 \in 
\mathcal{Q}_{\sss N} (| \mu^1 |)
\\
\mathbf{v}^2 \in 
\mathcal{Q}_{\sss N} (| \mu^2 |)
\\
\mathbf{v}^1 + \mathbf{v}^2 = 
\mathbf{v} }
$} &
*\txt{$
\mathcal{M}_{\sss N} (\mathbf{v}^1, \mu^1 )
\times 
\mathcal{M}_{\sss N} (\mathbf{v}^2, \mu^2 ) 
$}
\ar[d]^-{\theta_{\sss N}^k 
\times \theta_{\sss N}^k}
\ar[u]_-{\tau_{\sss N}}
& \\
*\txt{$\bigsqcup$ \\
$\substack{
\mathbf{v}^1 \in 
\mathcal{Q}_{\sss N} (| \mu^1 |)
\\
\mathbf{v}^2 \in 
\mathcal{Q}_{\sss N} (| \mu^2 |)
\\
\mathbf{v}^1 + \mathbf{v}^2 = 
\mathbf{v}
\\
r^1, \; r^2 \in \mathbb{Z}_{\geq 0}}
$} &
*\txt{$
\mathcal{M}_{\sss N}^k 
(\rho_{\sss N}^k (\mathbf{v}^1), 
\mu^1 , r^1)
\times
\mathcal{M}_2 
(\mathbf{v}^1_k ,
(\rho_{\sss N}^k (\mathbf{v}^1 ))_{k+1} ,
r^1)
\times$\\$
\mathcal{M}_{\sss N}^k 
(\rho_{\sss N}^k (\mathbf{v}^2) , 
\mu^2 , r^2)
\times
\mathcal{M}_2 
(\mathbf{v}^2_k ,
(\rho_{\sss N}^k (\mathbf{v}^2 ))_{k+1} ,
r^2)
$}
\ar[d]^-{(\Id \times \Id \times
\tau_2) \circ P_{23}}
& \\
*\txt{$\bigsqcup$ \\
$\substack{
\mathbf{u}^1 \in 
\mathcal{Q}_{\sss N}^k (| \mu^1 |)
\\
\mathbf{u}^2 \in 
\mathcal{Q}_{\sss N}^k (| \mu^2 |)
\\
\mathbf{u}^1 + \mathbf{u}^2 = 
\rho_{\sss N}^k (\mathbf{v})
\\
r^1, \; r^2 , \; r 
\in \mathbb{Z}_{\geq 0}}
$} &
*\txt{$
\mathcal{M}_{\sss N}^k 
(\mathbf{u}^1, \mu^1 , r^1)
\times
\mathcal{M}_{\sss N}^k 
(\mathbf{u}^2 , \mu^2 , r^2)
\times$\\$
\mathcal{S}_2
(( (\mathbf{u}^1_{k+1}, r^1 ),
(\mathbf{u}^2_{k+1}, r^2 )),
((\rho_{\sss N}^k 
(\mathbf{v}))_{k+1}, r ))
\times$\\$\times
\mathcal{M}_2 
(\mathbf{v}_k , 
(\rho_{\sss N}^k (\mathbf{v} ))_{k+1}, 
r)
$}
\ar[d]^-{\tau_N^k \times \Id}
& \\
*\txt{$\bigsqcup$ \\
$\substack{
\lambda \in 
\mathcal{Q}^+_{\sss N} 
(| \mu^1 |+| \mu^2 |)
\\
r \in \mathbb{Z}_{\geq 0}}
$} &
*\txt{$
\mathcal{S}_2
((\mu^1 , \mu^2 ), \lambda )
\times
\mathcal{M}_{\sss N}^k 
(\rho_{\sss N}^k (\mathbf{v}), \lambda , r)
\times$\\$
\mathcal{M}_2 
(\mathbf{v}_k , 
(\rho_{\sss N}^k (\mathbf{v}))_{k+1}, r)
$}
\ar@{<-} '[r] '[ruuuu] '[uuuu]_-{\Id \times 
\theta_{\sss N}^k}
& 
}}
\end{equation}
where $P_{23}$ is the permutation of the 
second and the third multiples in the 
direct product. It follows from 
constructions of the bijections 
$\tau$ and $\theta$
(or, more precisely, from comparison of
the right squares in the diagrams
\eqref{GLnGL2TTMMRestricted} and 
\eqref{GL2TTMMRestricted})
that the diagram \eqref{GLnGL2main}
is commutative.

The commutativity of the diagram \eqref{GLnGL2main}
together with Theorem \ref{GL2Theorem} implies 
that the bijection $\tau_{\sss N}$
commutes with the action of 
the Kashiwara's operators $\Tilde{e}_k$
and $\Tilde{f}_k$. Therefore
$\tau_{\sss N}$ is 
a morphism of crystals
\begin{equation}\nonumber
\tau_{\sss N} :
\mathcal{M}_{\sss N} (\mu^1) \otimes
\mathcal{M}_{\sss N} (\mu^2) 
\xrightarrow{\sim}
\bigsqcup_{\lambda \in 
\mathcal{Q}_{\sss N}^+ (| \mu^1 |+| \mu^2|)}
\mathcal{S}_2 ((\mu^1 , \mu^2 ), \lambda)
\otimes
\mathcal{M}_{\sss N} (\lambda) \; .
\end{equation}

This concludes the proof of  
Theorem  \ref{GLnTheorem}.

\subsection{Corollary}
\label{GLNCorollary}

The following proposition is a corollary
of Theorem \ref{GLnTheorem}.

\begin{proposition}
Let $\mu$, $\mu^1$, and 
$\mu^2 \in \mathcal{Q}_{\sss N}^+$.
Then
\begin{alphenum}
\item\label{GinzburgCorollary}
the $gl_{\sss N}$-crystal 
$\mathcal{M}_{\sss N} ( \mu )$
is isomorphic to
$\mathcal{L} ( \mu )$ (the crystal of the
canonical basis of the highest weight 
irreducible representation with highest
weight $\mu$);
\item\label{HallCorollary}
the cardinal of the set
$\mathcal{S}_2 ((\mu^1 , \mu^2 ), \lambda)$
is equal to the Littlewood-Richardson 
coefficient $c_{\mu^1 \mu^2}^{\lambda}$.
\end{alphenum}
\end{proposition} 
\begin{proof}
The proposition follows from
Theorem \ref{GLnTheorem},
Theorem \ref{JosephTheorem},
Proposition \ref{GLNHighestWeight},
and
Lemma \ref{GLNSumLemma}. 
\end{proof}
\begin{remark}
Note that \ref{HallCorollary} is the
Hall Theorem (cf .\ref{HallTheorem}).
Statement \ref{GinzburgCorollary}
could be deduced from the results of
Ginzburg concerning 
$gl_{\sss N}$-action in the top homology
of the variety 
$\mathfrak{M}_{\sss N} ( \mu )$ 
(cf. \cite{Ginzburg1991}).
\end{remark}

\subsection{Multiple tensor product and
Levi restriction}

There are two straightforward generalizations
of the constructions of this section.

First, one can consider more general tensor 
product variety $\mathfrak{T}_{\sss N} 
(\mu^1, \ldots , \mu^l)$, which
corresponds to a product of $l$
polynomial representations of $gl_{\sss N}$.
A point of this variety is a triple 
$(t, \mathbf{X}, \mathbf{F})$ consisting
of a nilpotent operator $t$, an $l$-step
partial flag $\mathbf{X}$, 
and an $N$-step partial flag $\mathbf{F}$, 
such that $t$ preserves both $\mathbf{X}$
and $\mathbf{F}$ and when restricted to
the subfactors of $\mathbf{X}$
(resp. $\mathbf{F}$) gives operators
with Jordan forms $\mu^1, \ldots , \mu^l$
(resp. $0$ operators). 
In this way one can prove a generalization 
of the Hall Theorem, saying that the number
of irreducible components of the 
Spaltenstein variety
$\mathfrak{S}_l ((\mu^1 , \ldots , \mu^l ), 
\lambda)$ is equal to the multiplicity of
the representation $L ( \lambda )$ in
$L ( \mu^1 ) \otimes \ldots \otimes
L ( \mu^l )$.

Second, one can generalize the notion of 
$gl_2$-restriction. Namely forgetting 
several subspaces in the flag $\mathbf{F}$
one can define restriction to a Levi 
subalgebra of a parabolic subalgebra of
$gl_{\sss N}$. In this way one would obtain
a bijection similar to $\tau_{\sss N}^k$
(cf. \eqref{GLnDefinitionOfTauNk}) 
relating tensor products for $gl_{\sss N}$
and for the Levi subalgebra.

\end{document}